\documentclass[a4paper,10pt]{article}
\usepackage{amsmath,amssymb}
\begin{document}

\title{SOME REMARKS ON VECTOR VALUED DISTRIBUTIONS}
\author{Tove Dahn (Lund University)}
\maketitle

\subsection{Invariance for change of local coordinates}

\subsubsection{Introduction}

In early global geometry, global means invariant for changes of local coordinates. We refer here to a global regular model, when the model has only regular approximations. Later, in global geometry, existence of regular approximations was discussed, we refer here to global very regular models.
It has been proved (\cite{Schwartz58}) that a global regular model is independent of choice of local coordinates, given a strong topology. We will discuss dependence of choice of coordinates in some other cases.

Assume $U \in \mathcal{G}$ a change of coordinates in a Lie group $\mathcal{G}$ and $T \in \mathcal{H}$, a space of distributions. We start with $<T,U_{j} I>=<U_{j}^{*} T, I>$, where we assume $U^{*} \mathcal{H} \subset \mathcal{H}$. Given $\Sigma f_{0} * \tau_{j} \delta_{0}$ represents a function in $L^{1}$
and when $<,>$ is given by convolution, the scalar product can be used to generate $\mathcal{H}=L^{1}$ through
$<U f,I>=\int f d U^{*}$ and $d U^{*} \bot I$ in $\infty$.
We assume a separation condition, for instance $\{ U \varphi  < \lambda \} \subset \subset \mathcal{B}$, for constants $\lambda$. A Schwartz type separation condition (\cite{Martineau}) gives that given 
$U$ locally reduced, $U^{*}$ is locally surjective. Given $d U \rightarrow d U^{\bot}$ is given by a subnuclear mapping, convexity is preserved. When we assume $E'$ is generated by $\mathcal{G}$, through a Stieltjes integral, we must assume $E'$ is nuclear and $<Uf,\phi>=<< f,dU>,\phi>$.

\subsubsection{Order of change of local coordinates}
Assume $d U$ BV (of bounded variation) on $\Omega_{\lambda}=\{ (x,y) \quad \mid d U \mid < \lambda \}$, a neighborhood of 0, for some constant $\lambda$ and consider the convex hull $T(\Omega_{\lambda})$. We assume there exists a restriction $d \tilde{U}$ with
$\tilde{\Omega}_{\lambda}=T(\{ (x,y) \quad \mid d \tilde{U} \mid < \lambda \}) \subset \subset \Omega_{\lambda}$. Thus we have on $\tilde{\Omega}_{\lambda}$, that $d \tilde{U}$ is invertible and Lebesgue is absolute continuous relative $d \tilde{U}$, that is $d \tilde{U}=d x /\alpha $, where $\alpha \in L^{1}(d \tilde{U})$. Schwartz condition means that $\tilde{\Omega}_{\lambda} \rightarrow T(\Omega_{\lambda})$ is compact.
Assume $U,U^{\bot} \in \mathcal{G}$, where $U^{\bot}$ is the result of $\psi : U \rightarrow {}^{t} U$ with respect to $<Uf,\widehat{g}>$, completed to $\mathcal{H}$. We will discuss $\psi \in \mathcal{L}_{c}(E,F)$, for locally convex and separated topological vector spaces E,F. The topology is given by uniform convergence over
convex, compact equilibrated sets (\cite{Schwartz57}).
Consider $W : \{ \phi \quad 0 < \mid U - U^{\bot} \mid(\phi) < \epsilon \}$, for a small positive $\epsilon$ and the convex hull $T(W) : \{ \phi \quad \mid U - U^{\bot} \mid(\phi) < \epsilon \}$. Assume existence of a regular approximation, for instance  $I \in \mathcal{L}_{c}(T(W),F)$, but $I \notin \mathcal{L}_{c}(W,F)$. Given $(U,U^{\bot})$ preserves reducedness, we have that $U^{\bot} \prec \prec U$, that is $U \neq U^{\bot}$ in $\infty$ and $I \notin \mathcal{L}_{c}$ in $\infty$. Note that, given $U-U^{\bot}$ analytic over $\phi$, we can consider $\phi(u,u^{\bot})$, as change of local coordinates. Example: given P corresponds to a partially hypoelliptic operator and $\tilde{\mathcal{L}}_{c}$ is the topology for convergence relative $d U$ BV over convex and radical domains. If further K corresponds to a parametrix to P, we have $K({}^{t} P \phi)=0$ is non-trivial, 
but $K_{N} ({}^{t} P^{N} \phi)=0$ is trivial, that is $\phi \in C^{\infty}$. Thus, $I \in \mathcal{\tilde{L}}_{BV,convex,equilibrated}$ and $I \notin \mathcal{\tilde{L}}_{BV,convex}$

\newtheorem{max_order}{maximal order}[section]
\begin{max_order}
Given $f,g \in \mathcal{H}$, consider a chain $L_{\mathcal{G}}(f,g)$ between f and g. The number of movements that define a chain is referred to as the order of the chain. For instance, consider for $G \in \mathcal{H}$,  $G(u_{1},\ldots,u_{k})$ over $L_{\mathcal{G}}$, then given $\frac{\delta G}{\delta u_{j}} \neq 0$ for some j, we have $d G \neq 0$. The chain has maximal order, k, if $\mathcal{H}$ can be generated by chains of order k. 
\end{max_order}

Assume A are chains in $\mathcal{G}$ maximal relative $\mathcal{H}$ and $(I)$, an ideal of functions regular and reduced over A. Assume in particular A algebraic over (I). Since $A \cap N(I)$ is algebraic, we have an analytic complement. We can assume for invariant sets in the complement, that they can be given by $d U_{j}=0$, $U_{j} \in \mathcal{G}$. Thus, for $\mbox{ rad }(I)$ we have partial algebraicity and the approximation property holds partially. Consider for instance $f(u_{1},\ldots,u_{k})$, with $d U_{k}=\alpha_{k} d U_{1}$ and $\alpha_{k} \in L^{1}(d U_{1})$, then existence of regular changes of coordinates is implied by the strict approximation property.
Note that precense of a change of coordinates through a maximal chain with invariant sets, is sufficient for f to be non-reduced. On the other hand, if $\mathcal{H}$ is generated by $\mathcal{G}$ over a base $f_{\alpha} \in L^{1}$ and $\mathcal{G}$ can be reduced to order 0, it is sufficient to consider translation.

Consider $U^{\bot}=(U-I) + V$ completed to $L^{1}$. We have that $U \in \mathcal{G}$ implies $U^{-1} \in \mathcal{G}$, but it is not necessarily reduced.  Assume $U^{\bot}=0$ implies $U=I -V$ with $R(V) \subset C^{\infty}$, then we have that U is topologically reduced in $\mathcal{D}^{' F}$, that is locally reduced modulo $C^{\infty}$.  

Givet two uniform spaces $X,Y$ (\cite{Bourbaki}), we have $Id : X \rightarrow Y$ uniformly continuous iff $X \subset Y$.
Consider $(U,U^{\bot})$, where we assume for some open set $\Omega$, that neighborhoods of $R(U) \cap \Omega$, are formed by movements of the same character as $U,U^{\bot}$. \emph{Assume $X_{1}^{\bot}=R(U)^{\bot},R(U^{\bot})=X_{2}$ are uniform then $\psi : X_{1}^{\bot} \rightarrow X_{2}$ uniformly continuous iff $X_{1}^{\bot} \subset X_{2}$, that is $U^{\bot}$ surjective}. In particular, we have this if $R(U)$ has the approximation property and we consider analytic functions $\widehat{g}$ with compact sub-level surfaces $\bot R(U)$. 

Assume $\Delta(f)$ is the lineality, that is the set of translation invariant lines for f.
\emph{Invariant sets for f can be derived from invariant sets for $\mid f \mid$ using uniformities.}
\emph{Assume $\psi : \Delta(f) \rightarrow \Delta(\mid f \mid)$ with image a uniformity. Assume B the set where  $\Delta(f) \simeq \Delta(\mid f \mid)$, where B is the uniformity that is given by the inverse to $\psi$.} Through Hurwitz (\cite{Ahlfors}), given the mapping is analytic, we have $\Delta(f) = \Delta(\mid f \mid)$ over B. We note that in the plane $(\mid f \mid, \mid \widehat{f} \mid)$, then $\mathcal{G}$ can be given by at most 8 movements (\cite{Lie91}). 
Note that given $U - I$ analytic over $f$ and Uf absolute continuous with $(U f)' = V g'$, then we have $d U(f)=d I(f)$
iff $V g'=g'$.

Given $d U^{\bot} = -\rho d U$, with $\rho \in L^{1}(d U)$ and $U \rightarrow U^{\bot}$ subnuclear , we have that $dU^{\bot}$ is absolute continuous relative $dU$. When $F$ is ac and symmetric with respect to $U,U^{\bot}$, that is $F \sim \int d F(U,U^{\bot})$, we have that $d F=(1 - \rho) d_{U}F$, that is $d F \equiv 0$ over spirals $(F=const$).

Note that linear independence is not necessary for a global model, it is sufficient with a global pseudobase (\cite{Oka60}), such that for instance $\mbox{ ord} \{ f_{1}=f_{2} \}$ finite. If $d f_{2} = p_{1} d f_{1}$, given $p_{1} \in L^{1}$ and Schwartz type separation conditions, we have that $p_{1}=1$ is removable. Note that $d U_{j}=p_{j} d U_{1}$, why we do not have necessarily that $U_{j} \rightarrow {}^{t} U_{j}$ preserve order (cf. $<f,\widehat{g}>$).

Assume U,V linearly independent on a domain $\Omega$, where both movements are analytic. Assume $U f_{j} = V f_{j}$ with uniform convergence on compact sets, then we have through Hurwitz (\cite{Ahlfors}), $(U-V)f_{j} \equiv f_{j}(u) - f_{j}(v)$ and $f(u)=f(v)$ or $f(u) \neq f(v)$ locally in $(u,v)$. Given monodromy, in the first case, we can identify $u-v=0$ locally. More generally, if $U+V=I$ in $\mathcal{L}_{ac}$ (over absolute continuous functions) and $d (U+V)=(\alpha + \beta) d x$, $(\alpha + \beta) f \in L^{1}_{ac}$,
we can conclude that $U \simeq I-V$ in $\mathcal{L}_{ac}$.

\subsubsection{Regular continuation}
Consider $d U = \alpha d V$, with $\alpha \rightarrow 0$ in $\infty$, this defines a convergence domain over lines. 
Given U monotonous, we can consider $U_{ac} + U_{sng}$, where  we assume that $d U_{sng}=0$ defines a regular complement. Assume $U \mid T \mid \simeq \mid V T \mid$, where  $V,U$ have the same character . We then have, given U algebraic, that $U \log \mid T \mid^{\theta} \simeq \theta \log \mid V T \mid$, that is given the condition $U \sim V$ and the approximation property, the conclusion is that $L^{\theta}$ is generated by translation and rotation outside the polar. 

Assume $d U =\alpha d V$, such that $1/\alpha \rightarrow 0$ iff $R(U) \subset R(V)$.  According to Riesz-Thorin: (\cite{Malgrange59}) given $\frac{1}{p_{\theta}}=1- \frac{\theta}{2}$, we have that $\forall f \in L^{p_{\theta}}$, there is $g \in L^{1}$ and $h \in L^{2}$, with $f=g+h$, that is $L^{p_{\theta}} \subset L^{1} + L^{2}$. Assume $d (A+B)=(\alpha + \beta) d U_{1}$,
where  $1/(\alpha + \beta) \in L^{1}(d U_{1})$. When we assume $(A+B)^{-1} : L^{1} \rightarrow L^{p_{\theta}}$, we must have $1/(\alpha+\beta)^{p_{\theta}} \rightarrow 0$. 
Conversely, given that $(A+B) : L^{2} \rightarrow L^{1}$, we have that $(A+B)^{-1} g=(C+D) f=f \in L^{2}$, that is we can choose $C=I-D$. Assume $1/(\alpha + \beta)^{p_{\theta}} \rightarrow 0$ in $\infty$. Sufficient for this is that $\alpha=p/q$ and $\beta=q/p$ under the condition that $\alpha \rightarrow 0$ or $\beta \rightarrow 0$ in $\infty$, for instance that $p/q \in \mbox{ rad } (I_{red})$. In particular, consider $(\alpha + \beta)^{p_{\theta}} = \alpha^{p_{\theta}} + \beta^{p_{\theta}} + r$, where  r does not change sign ($\alpha \beta > 0$).

Assume $R(U)^{\bot} = \cup R(U_{j}^{\bot})$, with $R(U_{j}) \subset R(U)$, that is the leafs are defined by $U_{j}^{\bot}$. Assume further that the completion to $L^{\theta}$ preserves character, sufficient for this is analyticity and linear independence. Existence of $U_{j}^{\bot}$, can be motivated by annihilator theory. Considering iterations, we assume $U \overline{f} =V f$, for a conjugated V. In this manner $U^{\bot}$ is defined as multivalued.

 Concerning monotropy, assume $\log f \in L^{1}$, such that we have existence $f_{x} \in C^{\infty}$ with $\mid f_{x} -f \mid < \epsilon$.
Then U is algebraic over $f_{x}$. Further, given $d \log f = \alpha \in C^{\infty}$, then we have that given $d^{2} f=0$ and $\alpha'=0$ on closed curves $\sim 0$, that
$\alpha d f$ defines a closed form. Assume $d U^{\bot}=\gamma d U$
with $\gamma \in \mathcal{D}_{L^{1}}(dU)$, given $\gamma'=0$ denotes max points on R(d U), we have that
$\mid dU^{\bot} \mid \leq C \mid d U \mid$ on $R(d U)$. Given the property of approximation through truncation, $\gamma$ can be used as a continuous deformation, that is $U^{\bot} \sim U$. Existence of $U$ such that $\forall f,g \in \mathcal{H}$, $Uf =g$  is  a lifting principle, that is $L^{1}$
has a lifting principle relatively $U_{1},U_{2}$. Assume $\{ d U \} \simeq \widehat{E}' \simeq \mathcal{G}$, where  $\widehat{E}$ is  a domain for U analytic in $\mathcal{G}$, that is ${}^{\circ} \widehat{E}'$. Given $\widehat{E}$ is Hausdorff, we have $\overline{L^{p}}(E) \simeq L^{p}(\widehat{E})$ (\cite{Schwartz57}). We can consider the continuation to $F$ such that U is algebraic over F, in particular U I = I U, that is $P_{1} \leq U \leq P_{2}$ locally over $F$.

\subsubsection{Conjugated invariant sets}

Assume F symmetric with respect to the reflection axes $L_{1},L_{2}$ and that $L_{1} \rightarrow L_{2}$ is bijective, then we have that $F$ is symmetric.
More precisely $(x,y) \underrightarrow{L_{1}} (x_{1},y_{1}) \underrightarrow{L_{2}} (y,x)$, for instance $z \rightarrow \overline{z} \rightarrow \overline{z}^{*}$,
where  F is separately symmetric, that is $F(z)=F(\overline{z})$ and $=F(z^{*})$. In the case where we do not have a bijection $L_{1} \rightarrow L_{2}$, for instance $z \rightarrow 1/z$ is not bijective, we use hypo continuous F. Under a hypo continuous mapping, the image of products of equi continuous sets, is relatively compact (\cite{Schwartz57}).
Consider $\sigma-\mu$ decomposable sets. For instance $d \mu (L_{2}) = \alpha(L_{1} \rightarrow L_{2}) d \sigma(L_{1})$ (\cite{Lie91} Lie factorization over conjugated transformations). In particular $L_{2}=L_{1}^{*}$, that is $d \mu( x^{*})=\alpha(x \rightarrow x^{*}) d \sigma (x)$ (\cite{Lie96} contact transforms, Lie). 
Note that given $I_{T-I}$ nuclear, we have that $I_{T} - I=I_{T-I}$.
If also I is  nuclear , then  $I_{T}$ is nuclear. 

\newtheorem{konjugerad}{Conjugated semi-norm}[section]
\begin{konjugerad}
Assume p,q semi-norms, with $q \leq p$. Assume existence of a semi-norm h, with $h^{2}=0$
algebraic and $\mid p(x)- h(1/x) \mid < \epsilon$ in $\infty$, (cf p preserves constant value in $\infty$).  Then we have given $U^{\bot}=(I-U) + V$, that $q(U^{\bot}) \leq \Sigma h(U^{j}) + p(V) + \epsilon$, as long as $h(U^{j}) \rightarrow 0$, when $j \uparrow \infty$ and we have $\mid q(U^{\bot}) - h((I-U)^{-1}) \mid < \epsilon$, given  $p(V)=0$.
\end{konjugerad}

When $U \rightarrow U^{\bot}$ is considered as conjugation, we assume $d U^{\bot} / d U \neq 0$. Assume further, for $T \in \mathcal{E}'$, that $\Delta$ is given in $\mathcal{G}$, with $<T, \varphi> \simeq <\Delta T,\psi>$, for $\varphi,\psi \in \mathcal{B}$ (\cite{Schwartz57}). Given $\Delta$ subnuclear,  we have ${}^{t} \Delta \psi \in \mathcal{B}$, given $\Delta$ nuclear, we have ${}^{t} \Delta \psi \in \dot{\mathcal{B}}$. Thus we can chose T such that $<T,\varphi>=0$, in $\mathcal{D}_{L^{1}}'$. Note that according to Lie, relative conjugation
it is sufficient to consider $<> \sim const$ (\cite{jag19}), for the orthogonal to be defined.

\emph{We can define an annihilator theory relative first surfaces: $X^{\circ}=\{ S \quad (\psi,S)=1 \quad S \in X' \quad \psi \in X \}$ and
${}^{\circ} X=\{ \psi \quad (\psi,S)=1 \quad S \in X' \quad \psi \in X \}$ (\cite{Nishino68}), given isolated singularities, we can prove that closedness for X, means  ${}^{\circ} X^{\circ} = X$}, for instance $X=L^{1},H$. Given S is scalarly absolute continuous, $d (S,\psi)=0$ implies $S \in X^{\circ}$. 
Finally, consider $X=\{ \psi \quad d \mid U \mid (\psi)=0 \}$, that is a domain over which $\mid U \mid$ is analytic. Given $\mid U \mid$ is absolute continuous (that is isolated singularities), we get the annihilator theory as above. 

\subsubsection{Integral representation}

Assume $\mathcal{H}$ normal, nuclear with $\gamma$-topology (\cite{Schwartz58}, prop 4, pg 41) and the strict approximation property, assume the strong dual is nuclear and quasi-complete. We can then  define a scalar product, hypocontinuous with respect to bounded sets in $\mathcal{H}(E)$ and compact sub sets in $\mathcal{H}'(F)$.

\newtheorem{Nevanlinna}[max_order]{Global regular model}
\begin{Nevanlinna}
Under the conditions in prop. 4 (\cite{Schwartz58}), we assume $<T,e'>(\varphi)=<T(\varphi),e'>$, that is $\mathcal{H}$ is in the algebraic dual to $E'$ and $E'$ is generated by $\mathcal{G}$. As long as $UT \in \mathcal{H}$, when $T \in \mathcal{H}$, we have that when $\mathcal{H}$ has the property of approximation by regularization and truncation, the representation is independent of choice of local coordinates in $\mathcal{G}$
\end{Nevanlinna}
(\cite{Schwartz58}, prop. 5, pg. 43).
Assume $E' \simeq \mathcal{G}$ and that for instance the mapping $U \rightarrow U^{\bot}$, maps $E' \rightarrow F'$. Given U reflexive, but not projective, we have for $U^{\bot}=(I-U) + V$, that $V^{\bot}=V$. For instance, starting from $\mathcal{D} \rightarrow E$ finite dimensional, we can choose  $E=F^{\bot}$. When $F=R(A)$, we have that E is  polar to A. Given $V,V^{\bot}$ are analytic, we have $\Gamma=\{ V=V^{\bot} \}=\{ 0 \}$. If they are only linearly independent continuous movements, we can have segments in $\Gamma$. In particular given $(V - V^{\bot})f < \lambda$, then $\mid x \mid \rightarrow \infty$.

\emph{Schwartz criterion for $\mathcal{D}_{L^{1}}'$, can be generalized to $\mathcal{D}_{L^{1}, \mathcal{G}}'$. For instance given $T \in L_{loc}^{1}$, we can give
$<T,d U> \in \mathcal{D}_{L^{1},\mathcal{G}}'$, as long as the coefficients to $dU$ preserves integrability.} Alternatively, given $L^{\theta}$ is generated by 
$U \in \mathcal{G}$, we can define $\mathcal{D}_{L^{\theta}}'$.

Given $E'$ is defined by derivatives and $\mathcal{H}=L^{1}$, we have $\mathcal{H}(E) \subset \mathcal{D}_{L^{1}}'$ or given $T \in L^{1}$, we have $<T,e'> \in \mathcal{D}_{L^{1}}'$, in this case we assume $\phi \in \dot{B}$. $\forall H \in \mathcal{D}_{L^{1}}'$ we have existence of a vector $T_{\alpha} \in L^{1}$, such that $<T_{\alpha},e'>=H$. Given $\phi \in L \subset \dot{B}$, we have existence of $H \in \mathcal{D}_{L^{1}}'$, such that $H(\phi)=0$, when L is closed we have that
$L^{\circ}=\{ H \in \mathcal{D}_{L^{1}}' \quad H(L)=0 \}$ and given $\phi \in \dot{B}$, we have $\phi \in L$.

Given $<f_{T},\widehat{g}>$ absolute continuous, we have $\int d <f_{T},\widehat{g}> \sim <f, \widehat{g}>$, for instance $\frac{d}{d T} <f_{T},\widehat{g}>=\int \alpha_{T}' < f, \widehat{g}> d T$, with $\int \alpha_{T}' d T=1$. 
Assume $T^{\bot}=(I-T) + V$, with $V \neq 0$ pseudo-local, then we have $<\varphi,T + T^{\bot}>= <\varphi,V>+<\varphi,I>$. Assume $V(\varphi)=0$ implies $\varphi \in C^{\infty}$, then T is projective modulo $C^{\infty}$. When $T^{-1}$ is very regular, we have $(I - T^{-1}) \in C^{\infty}$ and given $T : C^{\infty} \rightarrow C^{\infty}$, we have that T - I is regularizing, that is $T$ is very regular (in $\mathcal{D}^{'F}$).

Given $E^{j}$ is very regular $\forall j$, we have $(I-E)^{-1}$ is very regular  (a fixed singularity). 
Given  $TE = I$ (modulo $C^{\infty}$), then TE can modulo $C^{\infty}$ be considered as of type 0, symbols that constitute a convolution algebra. Given $TE = I$ (modulo $C^{\infty}$) implies $E=I$ (modulo  $C^{\infty}$), we have T has type 0. 

 Assume G a fundamental solution relative $I=I_{red}$ with $R(I_{red})^{\bot} =R(V) \subset C^{\infty}$, then G is parametrix outside the polar. Given further $V(\phi)=0$ implies $\phi \in C^{\infty}$, then a very regular  G can be continued to very regular in $\mathcal{H}(E)$. 

A disk-neighborhood can be generated by $(U_{1},U_{2})$. Assume $d U=\alpha d U_{1}$,$d U^{\bot}=\beta d U_{2}$. The domain for $(U,U^{\bot})$ is assumed to satisfy the condition that  $(\alpha,\beta) \in L^{1}(dU_{1},dU_{2})$. Note that
if we restrict U to acting on a disk, we have $U I = I U$, implies $P_{1} \subset U \subset P_{2}$, that is $\{ U = d U = 0 \}$ has Lebesgue measure zero. Note also $R(U)=\{ U f \quad f \in \mathcal{H} \} \simeq \mathcal{H}(E)$.

\emph{Assume $L(f,g)$ a connected path between f and g, then we have, given $f,g \bot \psi$ with $\psi$ in a nuclear space, that 
$<L(f,g),\psi>=0$}. \emph{In particular, given $R(U)^{\bot} = R(V)$, for some movement V, we have that V can be chosen as subnuclear}.  

Assume $L(\varphi,\psi)$ path between $\varphi,\psi$ and $\Delta : \mathcal{H} \rightarrow \mathcal{H}^{\bot}$ subnuclear , that is preserves convexity, then there is a corresponding path in the orthogonal.
More precisely, assume $<U \varphi,\phi_{1}>=<U \varphi,\phi_{2}>=0$ implies existence of path $U^{\bot}$ between $\phi_{1}$ and $\phi_{2}$, under the conditions that $\phi_{1},\phi_{2} \in R(U)^{\bot}$ and $U \rightarrow {}^{t} U$ continuous. Assume $d {}^{t} U=\alpha d I$. The condition $\alpha > 0$ locally, means that the movement does not change orientation locally. The condition $\alpha' > 0$ implies in particular that $U$ is absolute continuous locally. 

Assume $B(f,f')$ bilinear , and that a regular chain L between f,g defines a leaf for B. Assume $L(f',g') \bot L(f,g)$ with respect to B, and further that $(f',g') \sim 0$, that is $L$ gives a projective decomposition that can be reduced continuously to a trivial polar, that is we have ``surjectivity''. Given $(f,g) \sim (f',g')$ through conjugation (duality), we have that every chain between f,g has a corresponding chain between $f',g'$. Given $(f',g') \in R(V) \subset C^{\infty}$ we have a regular  chain. 
Given $U^{\bot}=(I- U) + V$ with $R(V) \subset C^{\infty}$ and U pseudo local, we have that
$(I - U)^{-1} U^{\bot} - I \subset C^{\infty}$, where  $N(U^{\bot})$ corresponds to a 1-polar. Note that 
$(I - U)^{-1} (I - U)f =0$ means that  $f(0) =  0$.

\subsection{Maximum principle}
\subsubsection{Introduction}
The criterion we use for inclusion $\mathcal{H}_{1}(E) \subset \mathcal{H}_{2}(E)$, related to weighted $L^{p}$-spaces, is that the quotient of weights tend towards 0 in $\infty$.

Starting from $<T(\varphi),e'>=<T,e'>(\varphi)$, where  $T \in \mathcal{H} \subset E^{' *}$, we see that $E_{1}^{'*} \subset E_{2}^{'*}$ induces the inclusion $E_{1}' \supset E_{2}'$, for instance $d U_{1} = \alpha d U_{2}$, with $\alpha \in \dot{B}$. But since $dU_{1}=0$ does not imply $d U_{2}=0$, the order of $E_{2}' \geq$ the order of $E_{1}'$.

Given Q a semi-norm on a separated space E, we can we define $N_{Q,p}(f)=(\int_{X}^{*} (Q(f(x))^{p} d \mu(x)))^{1/p}$ (\cite{Schwartz58}) and $\Delta^{p}(E)$ functions finite with respect to N. The closure of $\Delta^{p}(E) \cap \mathcal{E}'$ is not separated. The associated separated space is denoted 
$\overline{L^{p}}(E)$. We have $\overline{L^{p}}(E) \subset \overline{L^{q}}(E)$, when $p \leq q$. Given p,q semi-norms, corresponding to $dU_{1}$ and $d U_{2}$, for instance $p(f)=\mid \int f d U_{1} \mid= q({}^{t} \alpha f)$, where  $\alpha \rightarrow 0$ in $\infty$, we assume that $p/q \rightarrow 0$ in $\infty$. 

Assume $\mid < p(f), d U > \mid  \simeq \mid < f , d \mid V \mid > \mid$, for a semi-norm p. Thus when U is analytic it can be related to an harmonic measure. For example, given 
$\frac{\delta f}{\delta x} = \beta_{1} \frac{\delta p(f)}{d x}$,$\frac{\delta f}{\delta y} = \beta_{2} \frac{\delta p(f)}{\delta y}$ are well defined. Assume the coefficients to d U are $(\xi,\eta)$ and the ones to 
$d \mid V \mid$ are  $(\xi_{1},\eta_{1})$, so that the movements can be related through $\frac{\xi_{1} \beta_{1}}{\eta_{1} \beta_{2}} = \frac{\xi}{\eta}$. Further, we can construct a  semi-norm q, such that $\int p(f) d U_{T} \simeq \int q (\frac{d}{d T}(U_{T} f)) d T$, for a parameter T, given $q(U_{T} f)=0$ on boundary to the domain. According to the above, $q^{r}() \in L^{1}$ implies $q() \in L^{r}$, but $q(d U_{T}(f))=0$ does not imply $U_{T}$ analytic over f. Given a completing domain, where the norm is given by q, $U_{T}$ can be seen to be analytic under these conditions. Under the condition on the mean, $M(q^{r}(\frac{d}{d T} (U_{T} f))) \leq q^{r}(U_{T} f)$, we have that convergence for $U_{T}f$, with respect to $q^{r}$, implies convergence in $\overline{L^{r}(E)}$.

\subsubsection{Maximum principle}
The concept of dimension assumes linear independence for base and precense of 0, since $0 \notin C$, cylinder web, we prefer to consider order. In particular, for $f,g \in \mathcal{H}$ and $U \in \mathcal{G}$, we consider regular chains on the form $U f=\Sigma U_{j} f_{j}=g$, where $f_{j} \in \mathcal{H}$ constitutes a base for the change of coordinates. The movement is regular, if it is regular over base elements.

\newtheorem{max_princip}[max_order]{Maximum princip}
\begin{max_princip}
It is  to determine hypoellipticity, sufficient to consider chains of maximal order.
More precisely, if $F'$ are chains of maximal order between f,g that are given by U and if $E'$ are shorter chains and if E has the approximation property, given $I \in \overline{\mathcal{L}_{c}}(E,F)$, then F has the approximation property, that is U is algebraic.
\end{max_princip}
(\cite{Schwartz57}, Proposition 2, pg. 7, the closure is in $\mathcal{L}$ with topology induced by $\mathcal{L}_{c}$).
Note that U is locally such that $d U \sim \Sigma_{j \in I} dU_{j} \sim \Sigma_{j \in I} \alpha_{j} d U_{1}$, 
for some finite set of indexes and where $d U_{j}$ can have higher distributional order..
A necessary condition for $(U,U^{\bot})$ to preserve hypoellipticity is that $U^{\bot} \prec \prec U$ (\cite{jag20}). We assume also that the real part of the symbol does not change sign on a connected set $\ni \infty$.

Assume $T \in \mathcal{H}$. We denote with $\mbox{ ord } < T(x),e'>$, the number of movements in $E'$, that are necessary to generate $\mathcal{H}(E)$ locally. We assume $T*\varphi \in L^{1}$ implies $\mbox{ ord } \mathcal{H}(E)=0$, that is, is generated by analytic poly cylinders, on the support of T. Further, $T * \varphi \in L^{2}$ implies $\mbox{ ord } \mathcal{H}(E) > 0$.  Note that over $\mathcal{D}_{L^{2}}$, we have that U acts algebraically, that is assume $U \mid T \mid = \mid V T \mid$ with $\mid T \mid \in \mbox{ rad }(I)_{L^{1}}$, then we have that U can be limited to translations,
that is $\mid V T \mid^{N} \in L^{1}$ is of order 0 outside the zero-space.

Further, since $\mid \mid VT \mid - \mid T \mid \leq \mid V T -T \mid$, we have that 
$U \mid T \mid^{2}=\mid T \mid^{2}$ does not imply $V T = T$, that is $\mid T \mid^{2}$ reduced in $L^{1}$ does not imply $T$ reduced in $L^{2}$. Since $L^{2}$ is not nuclear, V can not be determined from U uniquely, and given a suitable separation condition, we do not have that V is surjective in $L^{2}$.

Assume $\phi$ real and does not change sign on connected components $\ni \infty$. Assume $\Omega= \{ \phi >0 \} \sim \{ F < \lambda \}$ is semi-algebraic $\ni \infty$. 
Then $\Omega \ni z \rightarrow 1/z \in \Omega$ continuous, does not preserve compactness for $\Omega$. 
Given $F(z) < \lambda$ implies $\mid z \mid < C$, we do not simultaneously have that $\mid 1/z \mid < C$. 
\emph{
But, given $\mid z \mid + \mid 1/z \mid \leq \mid  F(z,1/z) \mid$, we have that $F$ is  downward bounded separately in $z,1/z$. Thus F is ``hypo-reduced''.}  Assume further $F=\widehat{G}$, then we have
that the inequality implies $(1 + \mid z \mid^{2}) \leq \mid z F \mid$, thus given a condition on existence of j such that $\delta^{j} G$ is hypoelliptic (cf. very regular boundary),
we have that $F_{j}(z,1/z)$ is hypo-reduced as above.
Consider $\phi(z)=\mid z \mid + \mid 1/z \mid$, as long as $\mid z \mid$ is a monotonous deformation, $\phi$ has relatively compact sub level surfaces. 

When $L_{\chi}(\varphi)=<\chi,\varphi> \sim Id(\varphi)$ we have that $L_{\varphi} :\chi \rightarrow \check{\chi}$, that is the symmetry operator (\cite{Schwartz57}). Note that Id is optimally non-reduced, that is corresponds to a symmetric polar.  On D we have that $\varphi \equiv 0$ implies $\varphi^{*}$ (an algebraic dual) symmetric on $D^{' \circ}$ $(\varphi \rightarrow \varphi^{*}$). Note that determination of the spiral clustersets when $D^{' \circ}$ is the polar, is  in some cases unsolved. An algebraic polar implies $\varphi$ not identically zero on a disk (a necessary condition, given $\varphi$ is bounded and analytic, is that $\varphi$ is not identically zero on a line).

Given $f'$ convex, we have $\int f' d I = \int f'' d x = f'$, that is $I=I_{red}$. When $f''=\alpha f'$, we have $\alpha d x=d I$. 
Given $f \in X$, are convex functions, we have $\int f' d x=f$, that is I=$I_{RED}$, as opposed to evaluation in a point. $\delta_{x}(\varphi)=0$ if $\varphi(x)=0$, but $\delta_{x}$ has a translation invariant symbol, that is we have precense of a polar.
Let $\tilde{X}$ be a a convex continuation of X, then we have $\int_{\tilde{X}} f d I=\int_{X} \tilde{f}d I$, for instance $t f + (1-t)\int f' d I$. Thus I(f) can be
deformed continuously to f, analogous with an approximative identity.
On $T(X) \supset X$, for example $\int_{T(X)} f d I = \int_{X} f' d x + \int_{\Gamma} f d x$, 
that is given f=0 over $\Gamma$, we have I f = f. 
Consider $I$ completed to $I-\gamma : \mathcal{D}' \rightarrow \mathcal{D}_{L^{1}}'$.
 then for suitable $\gamma$, we have convex ($L^{1}_{ac}$, ramified) neighborhoods.

\subsection{Invertible movements}

\subsubsection{Resolution of identity}
Resolution of identity in this article, assumes a regular  covering.
Assume $S_{j}$ planes complementary to $R(U)$. Assume $S_{j}^{\bot} \bot$ $S_{j}$. Assume $S_{j}^{\bot}$ intersects  a tangent dU. The system R(U) does not necessarily have the property of approximation, but R(d U) can be completed to a system with this property.
More precisely, $\{ d U \} \cap S_{j}^{\bot} \neq \emptyset$
that is we have continuation of tangents $d U$ to $S_{j}^{\bot}$, sufficient for this is that $R(d U)^{\bot}$ has a projective decomposition, that is $S_{j} \cap S_{j}^{\bot}=\{ p \}$, p point. When R(U) is given by (integration of) $R(\Sigma d U_{j})$, we can get an approximation principle for the completed system.

Consider $\lim_{y \rightarrow x} \int \varphi(x,y) d y$, for x fixed on $S_{j}^{\bot}$ and $y \in S_{j}$. Given a normal tube (regular  covering), we have that the limit is not dependent of choice of starting point. Assume R(dU) completing and $d U \cap S_{j}^{\bot}$ a closed segment and $d U(0)=0$, then we have $E_{TR(dU)}$ is a Banach-space. Assume $S_{j}^{\bot} + S_{j}=I$ absolute continuous, then the condition for projectivity (local) is $d S_{j} \bigoplus d S_{j}^{\bot}=0$.

According to Lie, chapter 9, sats 8, given three equations represent a combination of movements,
we have that the corresponding integrals are  linearly independent. Assume the equations du=0, dv=0 and $d v - \Phi(u,v) d u=0$. Through a change of variables, when $\Phi$ is integral,  $\frac{\delta \Phi}{\delta u}=-\frac{\delta \Phi}{\delta v}$, and can we identify $d \Phi(u,v)$ and $d \Phi(u-v)$.
Given for instance $\Phi$ absolute continuous in u,v, we have that $\Phi(u,v)$ can be identified with $\Phi(u-v)$. 

The fundamental formula (\cite{Schwartz58}): $\varphi ._{\iota,\Delta} T = \varphi ._{\iota,I} S_{A}$, where I is  the identity, $\varphi \in \mathcal{D}_{A}$ and
$S_{A} \in (\mathcal{D}_{A})_{c}'(F,\beta_{0})$, where  $\beta_{0}$ are bounded sets, completing for F. Given $\Delta T=S$ and S invertible, we have that T is  scalarly invertible.
Note that $U \Delta T = \Delta U T$ on a closed segment, does not exclude precense of 
$R(U)^{\bot} \cap R(U^{\bot})$.  

\subsubsection{Invertibility}

Under the condition $T \rightarrow 0$ in $\mathcal{L}_{c}(\mathcal{H},\mathcal{H})$, we have $<T,\varphi> \rightarrow 0$ on convex, equilibrated, compact sets. In the same manner, $T \rightarrow I$, given $<(T-I),\varphi> \rightarrow 0$ in the same topology. Note that this does not imply that T is invertible, that is $T=I$ does not imply $\varphi=0$. The condition $T(\varphi) \rightarrow 0$ and further $T \rightarrow I$ implies $\varphi \rightarrow 0$, is  sufficient for T to have a continuous inverse.

Assume $\Sigma V^{j}=\Sigma_{N} V^{j} + A$, where  A is absolute continuous, that is $N(V^{j} - I) \subset N(V^{j+1} - I)$,
In particular, $\{ (V^{j} - I) < \lambda \} \subset \{ (V^{j+1} - I) < \lambda' \}$, for constants $\lambda,\lambda'$ and finite j. Assume $V^{j}(\phi) \sim V(\psi^{j})$, that is $V^{j}$ are absolute continuous over $(I)$ iff V is absolute continuous over $(I)_{j}$, we then have a domain for absolute continuity. Alternatively, assume $V^{j} - I \in C^{\infty}$ implies $V^{j+1} - I \in C^{\infty}$, that is $V^{j}$ are very regular in $\mathcal{D}_{L^{1}}'$.
 
Using for instance the spectral mapping theorem, given an algebraic change of local coordinates (over $\mathcal{D}_{x,y}$), invariant sets are preserved. 
Consider $(I)=(I_{0}) \bigoplus (I_{1}) \bigoplus (I_{2})$, where  $(I_{0})=$ $\{ \phi \quad U $ $\mbox{ absolute continuous over } \phi \}$, $(I_{1})=\{\phi \quad U^{j} \mbox{ absolute continuous over } \phi \}$, that is $\mbox{ rad }(I_{0})$ and $(I_{2})$ what is remaining. Then we have $(I-U)^{-1} \sim \Sigma_{N} + R$, where we assume R ac over (I). Obviously, $(I_{0}) \subset (I_{1})$. Since $(I-U)(\Sigma_{N} U^{j} + R)=I$, when we assume 
$d R(\phi)=0$ implies $U^{N+1}(\phi)=0$, we have $(I - U)^{-1} \sim \Sigma U^{j}$ a finite sum, that is, given $\phi \in (I_{1})$, with $d R(\phi)=0$, we have that $(I - U)$ is invertible.
Assume $U^{\bot}=(U-I) + V$, where $U,V \in \mathcal{G}$.

\newtheorem{rad_ac}[max_order]{Invertible radical}
\begin{rad_ac}
 Assume some iterate of U is absolute continuous and U is analytic over $\phi$ and that the polar is given by V, 
of the same strength as U,
then $U^{\bot}$ is invertible over $\phi$.
\end{rad_ac}

\emph{Consider $(I - \frac{U-I}{V})^{-1} \sim \Sigma  \frac{(U-I)^{j}}{V^{j}}$. 
Given $U^{j}$ is absolute continuous, when $j > N$ and $dU^{j}=0$, we then have $\frac{U^{j} -I}{V^{j}}=0$, when $V^{j} \neq 0$, that is it is necessary that V is non-trivial. Assume further the strength is
$U -I \sim V$, in the sense that $(U-I) / V^{2} \rightarrow 0$ in $\infty$, then we can derive the convergence
for $V^{-1}(I - \frac{U-I}{V})^{-1}$ from the conditions above.}

\subsubsection{Parametrices}

In $\mathcal{D}'(\mathbf{R}^{n})$ the polynomials with constant coefficients are surjective, that is given $U f \simeq P(D) f_{0}$ in $\mathcal{D}^{'F}$, with $f_{0}$ very regular,
we have existence (modulo $C^{\infty}$) of parametrices, that is $\exists V$ such that $VU \sim I$ in $\mathcal{D}^{' F}$. Given $P(D)$ hypoelliptic, the parametrices can be
chosen as very regular, that is modulo $C^{\infty}$, with trivial kernel. 
Consider $E : \mathcal{D}' \rightarrow \mathcal{D}^{' F}$, a parametrix to a polynomial operator. Given $E$ a projection operator (modulo $C^{\infty}$), that is $E^{2} \simeq E$, we can we assume $E^{\bot} \sim I - E$. Sufficient for this is that P(D) hypoelliptic. 

Assume $u : \mathcal{D} \rightarrow \mathcal{D}$ and simultaneously $u^{\bot} ; \mathcal{D}' \rightarrow \mathcal{D}^{'m}$, then we can have clustersets for $U$, even when ${U}^{\bot}$ has relatively compact sub level surfaces.  Given $U f=\int f d U$, we have $d U < \lambda d x$ implies $U f < \lambda$, when $f \in L^{1}$ (normalized), that is d U is  reduced, if U is  reduced over $L^{1}$. In the same manner for $\mathcal{D}_{L^{1}}$, even when U is algebraic over $\mathcal{D}_{L^{1}}$, we can have non removable zero-sets for ${}^{t} U$.
\newtheorem{phe}[max_order]{Partial hypoellipticity}
\begin{phe} 
Assume $f(u_{1},\ldots,u_{k})$ hypoelliptic over a non maximal chain. Assume $f \in \mbox{ rad }(I_{red})$ relative maximal chains. Assume  for a maximal chain, $U^{\bot}=(I-U) + V$, where  $I+V$
is algebraic, then $V$ can be represented as analytic, that is d V=0.
\end{phe}

Note that the condition that $f \in \mbox{ rad}(I_{red})$ relative maximal chains, implies that given an approximation property,
f is partially hypoelliptic, which in this context means that it is reduced relative translation. 
\subsubsection{Maximal rank}

Assume $PG=I_{red}$, for a hypoelliptic symbol P, and continue to $\tilde{I_{red}}$, with
corresponding $\tilde{G}$, parametrix to P. Then, $R(\tilde{G}) \subset \mathcal{D}^{' F}$, that 
is we have an approximation property. 

Starting from a covering $\{ \Omega_{j} \}$, where we assume $A_{j}f \in \mathcal{O}^{*}(\Omega_{j})$, that is analytic without invariant surfaces on $\Omega_{j}$ and f,
further $B_{i}A_{j}^{-1}f \in \mathcal{O}^{*}(\Omega_{j} \cap \Omega_{i})$ and $\Omega_{i} \cap \Omega_{j} \neq \emptyset$, then we have existence of a global U, such that $U A_{j}^{-1}f \in \mathcal{O}^{*}(\Omega_{j})$ $\forall j$. In particular, when $d U = \alpha d U_{1}$ and $\alpha \rightarrow 1$ regularly (without invariant sets), we have Oka's property (\cite{Range}). Over $\log f$, we have
for instance given $dI$ BV, $dI=dU + d V$, with $dU,dV$ BV and monotonous. Using Fourier duality, when $\Omega \ni 0$, we have a corresponding $\Omega^{*} \ni \infty$ and the condition $(U A_{j}^{-1})^{*}g \in \dot{B}(\Omega^{*})$. For analytic movements, we can assume $G(u_{1},\ldots,u_{k}) \in \dot{B}_{1} \times \ldots \times \dot{B}_{k}$ and we can determine a domain for ac and corresponding polar, for instance $I_{red} + V$ with $R(V) \subset C^{\infty}$.

\newtheorem{konvergens}[max_order]{Dependence of topology}
\begin{konvergens}
Assume X a domain for existence of a very regular fundamental solution G and consider $(U + W)$, with $U,W \in \mathcal{G}$ and $I \prec U$. Then we have that, $G {}^{t} (U + W)^{-1}$ corresponds to a very regular fundamental solution, only if $W \prec \prec U$ on X and $ W / U \in \mathcal{E}^{' (0)}(X)$.
\end{konvergens}

Given G-I of type $-\infty$, that is corresponding to regularizing action, we have $e^{c \mid x \mid} \mid G-I \mid < C$, for some positive constant c, when $\mid x \mid \rightarrow \infty$.
Thus the condition $e^{c \mid x \mid} \mid \frac{W}{U} G \mid < C$, when $\mid x \mid \rightarrow \infty$ is necessary and we do not necessarily have that G is very regular
under a change of coordinates, that generate the range. For parametrices E, given P is hypoelliptic, we have that, if $PE$ is of type $- \infty$ outside the diagonal, the same holds for E (Nullstellensatz). When E is Fredholm, we have that $\mbox{ ker} E$ is of type $-\infty$, and the condition is satisfied. 

\newtheorem{max_rank}[max_order]{Maximal rank}
\begin{max_rank}
 Assume P symbol to a d.o, hypoelliptic over chains of maximal order, assume the symbol G to a fundamental solution is in $\mathcal{H}(X)$, of order k. Assume (U+V) a change of variables of order k, then the resulting G is not necessarily very regular. Assume $X$ a poly cylinder in $L^{1}$, maximal for which G coresponds to a very regular distribution and that $I \in \overline{\mathcal{L}}_{c}(\tilde{X},X)$, where $\tilde{X}$ is a strict extension of $X$, then G is not very regular on $\tilde{X}$
\end{max_rank}

Assume $\mathcal{H} \ni f \rightarrow U f \in \mathcal{H}$ preserves order. When $(U+V)f=0$
with $V \prec \prec U$, then we must have $U \neq V$ in $\infty$, that is we do not have $U^{-1}V f=f$ in $\infty$. Presence of a non-trivial kernel $Vf=0$, that is continuation with 0, corresponds to presence of invariant sets (for $W=V+I$). Assume $R((U+V)^{-1})=D(U+V)$, that is $(U+V)$ preserves order, when $ord D(U+V)=ord R(U+V)$, which means absence of non-trivial kernels. In $L_{ac}^{1}$, this means $ord (U+V)^{-1}=0$ and when $\mathcal{H}$ has the approximation property, $(U+V)^{-1}$ can be reduced to order 0. Otherwise, we must consider completing sets.

Consider $PG=I$ in $\mathcal{L}_{ac}$, that is $(U+V) I (U + V)^{-1}=I$ in $L_{ac}^{1}$. Assuming the approximation property, the equation can be extended from $\mathcal{L}_{ac}$ to $\mathcal{L}_{c}$. Consider $d (U+V)=(\alpha + \beta) d U_{1}$ and a semi-norm q, such that $e^{C \mid x \mid} q(U_{1} G) \leq C'$, as $\mid x \mid \rightarrow \infty$ for some $C>0$, that is $U_{1} G$ has negative type and corresponds to regularizing action. For the extension to $\mathcal{L}_{c}$, we assume $q((U+V)^{-1} G) \leq e^{A \mid x \mid} q(G)$, as $\mid x \mid \rightarrow \infty$ for some $A>0$. Thus, $e^{(A+C) \mid x \mid} q(G) \leq C' e^{A \mid x \mid}$, as $\mid x \mid \rightarrow \infty$, is not necessarily of negative type and we do not necessarily have regularizing action. 

On the boundary for convergence we do not necessarily have twosided regularity, that is not regular surjectivity in $\mathcal{H}$ for $(U,V)$. Note that precense of essential singularities, means that hypo continuity is not sufficient for existence of twosided limits. 
Assume pluri complex convergence radius, that is L radius to A completing, $A \subset \tilde{A}$ and $\tilde{A} \cap L$ a closed segment, then we can assume that $(U+V)^{-1}$ is divergent on L, for instance $V/U < c$ on L. Here L corresponds to eigen vectors for $(U+V)$. In particular when $(U+V)^{-1} : \mathcal{D}^{' F} \rightarrow \mathcal{D}'$, then loss of derivatives and divergence is given for instance by $V/U \in \mathcal{E}^{' m}$. A discontinuous convergence is given by $(V/U)^{N} \in \mathcal{E}^{' (0)}$, when $N \geq N_{0}$, that is $V/U$  improves regularities. 

When X is a symbol to an operator that does not have maximal rank, there is a set for multivalentness in $X'$.
Relative norm, we assume $X^{\bot}=\{ x' \in X' \quad <x' , x>=0 \}$.   Given that dimension is  relative $C^{\infty}$, we consider $R(X)^{\bot}$ as generated by an ``algebraic'' functional, which is implied by the approximation property. Given $\{ F=c \}$ is given by a holomorphic function, there is a corresponding domain for multivalentness (\cite{Oka60}). 
Consider $f=e^{\phi}$, then we have $\phi \in C^{\infty}$ does not imply $f \in C^{\infty}$ and given $\phi$ has negative type, it does not imply f with negative type (regularizing action). Further, $\delta^{2} e^{\phi} = (\delta^{2} \phi + (\delta \phi)^{2}) e^{\phi}$ (cf. convexity outside the polar).
A sufficient condition, for $\phi$, complete with respect to maximal chains, to imply f complete with respect to maximal chains, is that $\mathcal{G}$ acts algebraically.
Note that $\{ f=const \}$ on sets $\{ \phi=0 \}$ and $\{ \phi = const \}$, but constant surfaces for $\phi$, do not have the same property as zero sets. However, given Schwartz separation condition and $\phi \in L^{1}$, the sets can be considered as equivalent. Note that invariant sets for $f$ corresponds to invariant sets for $\phi$,
that is $U e^{\phi} - e^{\phi}=0$ iff $e^{U \phi - \phi}=0$.

Starting with $\int U^{\bot} d f = \int f d U^{\bot} - \int_{\Gamma} U^{\bot}f$. Assume $\Gamma$ a line in the polar $U^{\bot} f=0$. Given $U^{\bot}$ is algebraic over the tangent space, we have outside the polar and over non trivial f, that $U^{\bot}f \neq 0$ a.e. Given $R(U^{\bot})$ convex $\subset H$ a Banach space, the polar defines a completing set.
Consider now $\{ U^{\bot} \varphi \} - \{ U \varphi \}^{\bot}$. 
Assume for instance X is defined by $\mathcal{G}$, such that $\{ U \varphi \}^{\bot} \simeq \{ V \varphi \}$, where  $V \in \mathcal{G}$, is not necessarily unique.
Further, $d U^{\bot} \sim p d x$ and $d U \sim q d x$, where  $p$ is reduced and $ p/q \rightarrow 0$ in $\infty$. Given $p / q$ algebraic, there is not precense of a spiral. $p \sim q$ implies possible precense of a spiral.  Assume  $H=R(U) \bigoplus R(U)^{\bot}$ and that U is monotonous
with $d U = \alpha d I$, $\alpha > 0$ locally, that is $\alpha_{ac} + \alpha_{s} > 0$ locally. Given $\int \alpha_{ac} d I(f) =I_{red}(f)$ locally, we have $U f=I_{red}(f) + U_{s} (f)$, with $d U_{s}(f)=0$ locally, note that $U_{s}$ is dependent of the topology. 
To conclude that $I \phi=0$ nbhd 0, implies $\phi \equiv 0$ ($\delta * \phi = 0$ implies $\phi(0)=0$), it is necessary that $\dim_{C^{\infty}} R(I)=0$.

\subsection{Topological convergence}
\subsubsection{Boundary condition}
Assume a very regular  boundary (\cite{jag13}), that is precense of regular approximations, in particular $\exists \varphi \in \mathcal{D}$, $T * \varphi \rightarrow T$ with $T \in \mathcal{H}$, where $\mathcal{H}$ is normal. The boundary can be defined through the condition $\{ F^{(j)} = const \}$ isolated points, for some analytic function $F$, which excludes precense of a spiral. Consider $f \in X / C^{\infty}$, that is
$f \sim \delta - \gamma$ with $\gamma \in \mathcal{D}$, where  $\gamma \rightarrow 0$, it is then sufficient to consider f of type 0 (real type). Given $f \in \mathcal{D}_{L^{1}}'$, we have that $f \sim \Sigma D^{\alpha} f_{\alpha}$, where  $f_{\alpha} \in L^{1}$ (\cite{Schwartz66}),
that is given $f_{\alpha}$ has distributional order  0, we have that $f$ has finite distributional order . Assume $\delta f=c f$ (that is $\delta \log f=c$) on segments,
but $\delta^{2 }f \sim c f$ implies isolated points, then we have $\delta^{2} \phi f \sim (1 + \delta \phi)(1 - \delta \phi) f$. Sufficient for this is twosided regular limits and we have a decomposition $(\delta \phi)^{2} \bigoplus \delta^{2} \phi=1$, where $f \neq 0$.

 Given $d v \in \mathcal{D}^{' m}$, we have existence of $d \mu \in \mathcal{D}^{' 0}$, such that  $d \mu = \alpha d v$, where  we choose $\alpha \sim 1/P$ (modulo $C^{\infty}$), for a polynomial P, such that $\alpha : \mathcal{D}' \rightarrow \mathcal{D}^{´  F}$ surjective, (\cite{Schwartz57}, the property ($\epsilon$)).  Alternatively, consider $\tilde{d v} \sim d v$ (modulo $C^{\infty}$), for instance $\alpha d v(\varphi) = 0$ implies $\varphi \in C^{\infty}$.
In particular, given $\alpha$ very regular, we can associate $d \mu$ to $d v \in \mathcal{D}^{' m}$.

Assume $bd A=\{ U = \lambda I \}$, that is given U is absolute continuous we can write $bd A = \{ d U=0 \}$. Given U analytic over bd A, the boundary can be given by a connected set. When U is analytic over A, we have $U \leq \lambda$, that is maximum is reached on the bd A.
Given $\Omega$ pseudo convex, we have that bd $\Omega$ is  of order  0 (cylindrical). Consider $T \in \mathcal{E}'(\Omega)$, such that $T \mid_{bd \Omega} \in \mathcal{E}^{' (0)}$ and $\Omega=\{ (x,y) \quad P(x,y) < \lambda \}$. For instance $S \circ (P,\overline{P})(x,y)=T(x,y)$, where  $P < \lambda$ iff $\overline{P} < \overline{\lambda}$ and $S \in \mathcal{E}^{' (0)}$. In this case we have $P=S^{-1} T$ with respect to an involution condition. 

Every differential operator on the real line, corresponds to a reduced measure, in particular when $\lambda \in \mathbf{R}$, we can assume $P(D_{\lambda}) f_{\lambda}=0$ implies $f \in C^{\infty}(\lambda)$, that is for finite values on $\lambda$, we can assume an analytic dependence of $\lambda$.

\subsubsection{Projectivity}

Given $U=I$, we have that $U^{\bot}$ is not necessarily uniquely determined. In the plane there are finitely many leafs.  In particular, given I is separately continuous on $(U,U^{\bot})$ and equi continuous on leafs of BV, then $I(U,U^{\bot})$ is considered as hypo continuous. Note that, given $(U,U^{\bot}) \simeq (U,\psi(U))$ gives a closed graph, under a Schwartz type separation condition, we have that $\psi$ is  continuous. Further, that given isolated singularties, we can we identify $\{ U=I \}$ with $\{ U=0 \}$, for instance given $U \phi \in L^{1}$, we have $\{ e^{U \phi} \sim I \}$ iff $\{ U \phi \sim 0 \}$. 

Consider $\psi : d U \rightarrow d U^{\bot}$, with $\tilde{\psi}$ very regular , that is $\tilde{\psi}(d U)=d U^{\bot} + d R$, with $d R$ regularizing, that is $d U - d U^{\bot}=d R \in C^{\infty}$ but $d R(\varphi)=0$ does not imply $\varphi \in C^{\infty}$. Given that R is very regular and can be chosen downward bounded, we have $(R - I)(\phi) \in C^{\infty}$ and $R(\phi) \in C^{\infty}$ implies $\phi \in C^{\infty}$. Note, given R is absolute continuous, $d R(\phi)=0$ 
(modulo $C^{\infty}$) implies $(R- I)(\phi) \in C^{\infty}$.

Given duality $L \rightarrow L'$ and $M \rightarrow M'$, we can continue $T : L' \rightarrow M \rightarrow M'$ to $T' : L \rightarrow M' \rightarrow M$, assuming that the mapping $M \rightarrow M'$ can be converted to $M' \rightarrow M$ (cf. two mirror model). In particular if M is closed, we have that ${}^{\circ} M^{\circ} \simeq M$. Assume $(L \otimes M)^{\circ}=\{ (L \otimes M)' \quad < l(\varphi),m(\psi) >=0 \}$, that is $l(\varphi) \bot m(\psi)$.
Thus $L^{\circ} \otimes M^{\circ} \subset (L \otimes M)^{\circ}$. Given $L',M'$ reduced, we have however ${}^{\circ} (L \otimes M)^{\circ}=\{ 0 \}$.
 Further, given $U T \in C^{\infty}$ and $U^{-1}$ pseudo local, we have $T \in C^{\infty}$, that is U is algebraic. 

Assume K very regular  with respect to $(U,U^{\bot})$ and a norm $\parallel \cdot \parallel$, then the condition $\parallel U - U^{\bot} \parallel= 0$, gives precense of a spiral, but
$\parallel U \parallel - \parallel U^{\bot} \parallel = 0$ implies Parseval (the movement is considered as an operator). Further, $\parallel d U - d U^{\bot} \parallel \rightarrow 0$, given U absolute continuous, means  $(U - U^{\bot})=I$, that is $U^{\bot}=I-U$ and we have projectivity.

Assume for spaces of distributions, $\mathcal{H},\mathcal{K}$, that $\Omega=$ nbhd 0 in $\mathcal{H}$ and that we have existence of $\Omega'=$ nbhd 0 in $\mathcal{K}$, with $\Omega' \subset \subset \Omega$ and $\psi : \Omega \rightarrow \Omega'$ (restriction),
for instance $L^{1} \rightarrow L^{1}_{ac}$. In particular $\psi(d U)=\alpha d U_{1}$, with $\alpha \in L^{1}(d U_{1})$ and $\alpha \rightarrow 0$ in $\infty$ or $d I' = \alpha d I$ with $\alpha \rightarrow 0$ in $\infty$. Then over a regular  chain in ($L^{1}_{ac})$, we have that I is reduced. Given instead $\alpha$ independent of some variables, the separations condition above is not satisfied. Assume conversely $\psi$ an extension and $\Omega',\Omega$ disks, with $\psi(\Omega') \subset \Omega$, then $\psi$ is sub-nuclear and for instance $\beta d I'=d I$, with $\beta$ downward bounded  (I is invertible relative the representation). Finally, consider $T(f)=\int f d U$, where  $\mathcal{H}(E)$ is generated by $E' \simeq d U,d V$. Given $T(f)=const$, when V is varying, we have for the corresponding regular chain, precense of a polar (multivalentness).

\subsubsection{Orthogonal decomposition}
When E is  quasi-complete, we have $\forall B$ bounded $\subset E$, that B is completing (\cite{Schwartz58}). 
For completing B, we have $F_{B}$ is a Banach-space, that is we can discuss orthogonal decompositions. 
Note, given the parallellogram  law $2 (U,U^{\bot})=\parallel U - U^{\bot} \parallel^{2}$, that is $(U,U^{\bot})(\phi)=0$ iff $U \phi=U^{\bot} \phi$. Note in particular when $U^{\bot}=(I-U) + V$, we have $\parallel I + V \parallel^{2} + \parallel I - V \parallel^{2}=\parallel U \parallel^{2} + \parallel U^{\bot} \parallel^{2}$. However an orthogonal decomposition is not necessary for a global model (pseudo-bases). 

 Consider $U^{\bot}(\varphi)=0$, with $\varphi_{1} \otimes \ldots \otimes \varphi_{n}$ dense in $\varphi$, for instance over a domain where  $U^{\bot}$ is analytic over $\varphi$ (pluri complex zero). Then $U^{\bot}$ can be represented
through $(V^{\bot}_{1},\ldots, V^{\bot}_{n})$, such that $U^{\bot} \varphi \sim (V^{\bot}_{1} \varphi_{1},\ldots, V^{\bot}_{n} \varphi_{n})$. In particular, given $\varphi_{1} \bot \varphi_{2}$, there is a decomposition through conjugated $V_{j}^{\bot}$, $(W,W^{\bot})$.

The condition of existence of $U \in \mathcal{G}$ analytic over a chain to f, means that
$\overline{\delta_{1}} \ldots \overline{\delta_{k}} f(u_{1},\ldots,u_{k})=0$. Note that maximal order  with respect to $\mathcal{G}$, does not assume that all elements are analytic (Oka's property). 

\newtheorem{ac_decomposition}[max_order]{Absolute continuous representation}
\begin{ac_decomposition}
 Assume U absolute continuous. Given $dU=0$ can be decomposed into linearly independent measures, $\cup \{ d U_{j}=0 \}$, $\cap N(d U_{j})=\{ 0 \}$ and given R(U) has the approximation property, then U has a projective decomposition over regular chains.
\end{ac_decomposition}

Assume $d U=\xi \frac{\delta f}{\delta x} + \eta \frac{\delta f}{\delta y}$ and $\frac{\delta G}{\delta y}=-\xi$ and $\frac{\delta G}{\delta x}=\eta$. Then we have 
$<d U(f),\phi>=<f, (\eta_{y} -\xi_{x}) \phi> + <f, d U(\phi)>$. Given $G$ regular , $\{ G \}$ locally 1-1 iff ${}^{t} \{ G \}$ surjective, gives the approximation property through regularization. 
Given $\Delta G=0$ we have that
$d U = d {}^{t} U$. In order to apply the partial Fourier transform, we must assume $\xi,\eta$ algebraic in x,y. In this case we have $(x,y) \rightarrow (\xi,\eta)$ subnuclear.
Given also $I$ is subnuclear, we have that the coefficients are in a nuclear space.
Assume, starting from $\{ f,G \}=0$, that $R(d f)$ is generated by $U \times U^{\bot}$, then $R(d G)$ is generated by ${}^{t}(U \times U^{\bot})$. Given reflexivity, ${}^{t} (U,U^{\bot}) \simeq (U^{\bot},U)$ that does not affect the order (dimension). 

\subsubsection{Polars}

Assume $\Sigma=\{ \alpha \beta=1 \}$, where  $d U=\alpha d I$ and $d W=\beta d U$, where $\alpha,\beta \in \dot{B}$. Assume L an oriented chain of maximal order, then we have given the approximation property, that $\Sigma$
gives invertibility.  Note that $d W(f) \in C^{\infty}$ implies $W f \in C^{\infty}$. 
Given $U,W$ pseudo local, then $U f \in C^{\infty}$ implies $W^{-1}$ pseudo local over $\Sigma_{\alpha \beta}$ and $W$ algebraic over $\Sigma_{\alpha \beta}$. 

Assume $A \subset B$ a disk neighborhood, where  B is defined through $\mid \xi \mid^{\sigma} d (\xi,\Gamma) < c$ (weighted disk neighborhood), where  we assume $\Gamma$ is defined as analytic. Then we have that $B \backslash \Gamma$ is analytic. Starting from hypoellipticity, we can define $\Gamma$ through chains of maximal order , that do not contribute micro locally, that is $\Sigma_{\alpha \beta \gamma}$ are trivial. Note for shorter chains, $\Sigma_{\alpha}$ trivial, does not imply $\Sigma_{\alpha \beta}$ trivial.

\newtheorem{ram_polar}[max_order]{Ramified polar set}
\begin{ram_polar}
Given $f/g=\alpha$ analytic, with $\alpha \rightarrow 0$ on a radius starting from the origin, then we have that $<\alpha g,g>_{L} \equiv 0$, that is $f \bot_{L} g$. Assume $\Omega$ is completing, then we can motivate $f \bot_{\Omega} g$. Assume $U^{\bot} \in \mbox{ rad } I(U)$, then we have $U^{\bot} \sim U^{N}$ (Nullstellensatz). The condition $U^{N} + U^{\bot} =I$ gives a ``rectifiable'' change of variables, or a ramified polar.
\end{ram_polar}

Using for instance Lindel\"of's theorem, given $\alpha \in \dot{B}_{u_{1},\ldots,u_{k}}$, we have that given $\alpha=0$ has a positive measure, on a line through 0, then $\alpha$ is zero on a disk. If we assume $\alpha=\alpha_{1} \dot \ldots \alpha_{k}$, we have motivated $f \bot_{\Omega} g$.
Given $\mathcal{H}$ a space of distributions, we have $\mathcal{L}_{c}(\mathcal{H}_{c}',E) \simeq \mathcal{H}(E)$. 
The topology for polars is given relative $\mathcal{L}_{c}$, that is uniform convergence on convex, compact and equilibrated sets. We can for instance assume on polars, that $f^{2} \in K$ does not imply $f \in K$. 
Or if T is in $\mathcal{L}_{c}$ on $0 < \mid U - U^{\bot} \mid < c$, but not in $\mathcal{L}_{c}$ on 
$\mid U - U^{\bot} \mid < c$, that is T is not uniformly continuous on the spiral, then we have that T is in the polar relative $\mathcal{L}_{c}$.

Assume $E'$ is generated by $\mathcal{G}$, $UT(\varphi)(\zeta)=T(\varphi)(\zeta_{T})$. Note that we do not need a global model for $\zeta_{T}$.
It is  sufficient that ${}^{t} U \varphi$ is surjective and that locally ${}^{t} U \varphi(\zeta) \rightarrow \zeta_{T}$ continuous. 

Assume $d V=\alpha d U$ continuous, where  $\alpha \rightarrow 0$ in $\infty$. Given $\varphi \in B$, we have $\varphi  \alpha \in \dot{B}$, why $d U \in \mathcal{D}_{L^{1}}'$ (with approximation property), that is given very regular boundary, where $\varphi \rightarrow 0$ in some direction, we have a decomposition of variables.
Assume $\phi \in (\dot{\mathcal{B}})'$ and completed to $\phi -I \in C^{\infty}$, 
then $\phi$ can be given by a hypoelliptic operator outside the polar.

\newtheorem{abel}[max_order]{Preservation of character }
\begin{abel}
 Assume $d U$ Stieltjes measures, defined according to Lie, with $C^{\infty} \supset R(U) \rightarrow R({}^{t} U) \subset C^{\infty}$ continuous. Assume $\Gamma=\{f  \quad d (U - {}^{t} U)(f)=0 \}$, that is the movement preserves character  over $\Gamma \subset \mathcal{H}$. Assume $\Gamma_{0}=\{ f \quad I \in \overline{\mathcal{L}}_{c}(U f    \rightarrow {}^{t} U f) \}$. then we have $\Gamma_{0} \subset \Gamma$.
\end{abel}
Consider the condition $\eta_{x} - \xi_{y}=0$ according to $\Delta G=0$. Given $U - {}^{t} U$  absolute continuous, we have that $\Gamma_{0} \neq \Gamma$ on non trivial sets.
Note I K = K I
implies $P_{1} \prec K \prec P_{2}$, for polynomials $P_{j}$. More precisely, assume K has symbol in H (holomorphic) and $K(e^{\phi}) \sim e^{K(\phi)}$, with respect to Exp-norm in the Fourier-Borel symbol space. Thus $K(\phi^{j}) \sim K(\phi)^{j}$, which means that the zero space to the symbol to K, can be given as geometrically equivalent with the zero space to polynomials. Thus, we have existence of $P_{1},P_{2}$, as above.

Consider $\Delta T=S$, for instance $\widehat{T}=S$ or $T^{\bot}=S$.
Starting from completing sets, with $F_{0}=\mathcal{G} E$, we can define $F_{0}^{\bot}$ through a Banach-norm. Assume $F_{1}=(I-U)E + VE$, where  $F_{1}=F_{2} + P$ with $P=R(V)$, then  we have $F_{0} + F_{2}=E$ and $F=E +P$. Note in particular, given $U \in \mathcal{G}$ then $U^{-1} \in \mathcal{G}$,
but when $U^{\bot}=(I-U) + V \in \mathcal{G}$, we do not have $\exists (I-U)^{-1} \in \mathcal{G}$. Consider $E' \rightarrow F$ and $F^{\circ}=\{ x' \quad <x,x'>=0 \quad x \in F \}$. Given F is closed, we have $\{ x \quad x' \in E' \quad <x',x>=0 \} \subset F$.

Consider the polar as an envelop, for instance $M \subset I(\otimes d \mu_{i})$, that is an ideal that is generated by tangents,
for instance that the polar is defined as a zero space to Lie-movements (analytic movements).
Note $<X(x,y),f(x)f(y)>$, means that given X symmetric, that $I_{X}(f_{x})=0$ iff $I_{X}(f_{y})=0$. Concerning linear independence, $<I,\varphi>=0$ iff $\varphi(0)=0$, that is we do not have $\varphi=0$ nbhd 0. Assume $\varphi \in Z$, then we have $T(\varphi)=0$ implies $T(\varphi) \in (E')^{\circ}$, implies $<T,e'> \in Z^{\circ}$. Note that given $T \in \mathcal{H}$, that is generated by $E'$ and
$<T,e'>(\varphi)=0$ implies $<T(\varphi),e'>=0$, then the implication $T(\varphi)=0$ means that we do not necessarily have a global regular model for $\mathcal{H}$.

\subsubsection{Decomposable sets}

Assume $B \subset L \otimes M$, where  $L=R(d U)$,$M=R(d U^{\bot})$. As long as  $I \neq B$, B can be represented through a projection operator.
Assume $B \subset L_{ac}^{1}(d \mu)$, Assume $d v=\alpha d \mu$, with $\alpha \in L^{1}(d \mu)$. We then have $\frac{d^{2} v}{d T^{2}} \geq 0$ iff $\frac{d}{d T} \log \alpha \geq - \frac{d}{d T} \log \frac{d \mu}{d T}$. Given $d v - 1 d x \sim d \mu$, means $d \mu= 1/(1 - \alpha) d x$, given $\mid \alpha \mid < 1$ is invertible.

Assume $\psi : d \mu \rightarrow d v$, for instance according to Radon-Nikodym (\cite{Riesz56}), $d v=\alpha d \mu$ with $\alpha \in L^{1}(d \mu)$. Assume $L \otimes M \subset L^{1}_{ac}(d \mu)$, for instance $d v$ absolute continuous with respect to $d \mu$ and v absolute continuous, then we have that the integral over $L \otimes M$ is defined.  We can assume on bounded sets, that $\psi$ is  defined modulo analytic action, that is $\psi U=U^{\bot} + V$, where  $d V=0$, thus modulo regularizing action.

Consider V : ST=1 and define $\tilde{V}=\{ y \quad T(y)=T^{-1}(x) \quad x \in V \}$.
Note that $S(x)=0$ implies $T(y)=0$, that is $x \notin V$, but $T(x)=1$ implies T(y)=1, that is we consider conjugated surfaces $x \sim y$,
that is $\exists x \rightarrow x_{0}$ with $T(x)=1$ and $\exists y \rightarrow x_{0}$ with $T(y) \neq 1$
(very regular boundary). 

\emph{The condition $T(x) T(y)=I$, can be written $e^{t + t^{-1}}=I$. Assume $t + t_{1}=I$, then we have $I=t_{1} - t^{-1}$. Given projectivity $R(t^{-1})^{\bot}=R(t_{1})$. The condition $N(t^{-1}) = \{ \phi \quad t^{-1}(\phi)=0 \} \subset C^{\infty}$ corresponds to a algebraic (removable) set. Given analyticity for T, we can consider
$x(T) y(T)=I$, that is $y(T)=1/x(T)$.} 

\subsubsection{Discontinuous movements}
 Given $\mathcal{H}$ scalarly of order m with property $(\epsilon)$ and $E$ 1-dimensional, then $\mathcal{H} \subset \mathcal{D}'$ is of finite order (\cite{Schwartz57}).
In the case of $\dim E > 1$, there are topologies involving a spiral, where the limit can be undetermined. I can be approximated by $I - \gamma$, of finite order, when $\gamma \in \mathcal{D}$ . Note that modulo $C^{\infty}$, operators of finite type can be considered as of type 0 (real type). 

When $\tilde{T}=T*\phi=0$ implies $\phi \bot R(T)$, we note that
given $\tilde{T} \rightarrow T$, we do not have $\delta \bot R(T)$. Sufficient for convergence is  $\tilde{T}$ downward bounded, as $\tilde{T} \rightarrow T$.
Assume for $\Omega$, $T(\Omega)$ the closed convex hull, that is the barrel to $\Omega$. Assume W a barrel $\subset \Omega$. Assume $p_{1}$ a semi-norm to $U \mid_{W}$,
$p_{2}$ to $U \mid_{\Omega}$ and $p_{3}$ to $U \mid_{T(\Omega)}$. A condition equivalent with inclusion (related to $L^{p}$) is  $p_{3} / p_{1} \rightarrow 0$ in $\infty$, for instance $p_{3} / p_{2} \rightarrow 0$, where $p_{2} / p_{1}$ is bounded. Given $p_{1},p_{3}$ algebraic (removable zero sets), we have hypo continuous convergence in $\mathcal{L}_{c}(\Omega,E)$.

Necessary for hypo continuity is  $I \circ e'=e' \circ I$ (\cite{Schwartz57}, Prop. 18), that is topological algebraicity.
Assume $F : E \rightarrow E$, according to $F ( \mid e \mid) \rightarrow \mid F(e) \mid$ equicontinuous. We have given $\mid F \mid (\mid x \mid, \mid y \mid) \in \mathcal{E}^{' (0)}$, that $F \in \mathcal{E}^{' (0)}$. For an equicontinuous mapping $E \ni (x,y) \rightarrow d U(x,y)$, we have that $\mid d U \mid<1$
corresponds to $(x,y) \in nbhd 0$, for instance $\{ (x,y) \quad \mid d U \mid < 1 \} \subset \subset E$ (\cite{Bourbaki}).
Assume existence of U, such that $T*U I \in \mathcal{E}^{'(0)}$, that has the approximation property. Further, $U_{1} I *T + U_{2} I * T \in \mathcal{E}^{' (0)}$ and given $U + U^{\bot}=I$ (outside the polar), we have $T \in \mathcal{E}^{' (0)}$. In particular, when $(U,U^{\bot}) f_{0}$ generates $L^{1}$ (outside the polar)
we have that ${}^{t} (U,U^{\bot})$ is  locally 1-1 (outside the polar).

\newtheorem{diskontinuitet}[max_order]{Discontinuous limits}
\begin{diskontinuitet}
Spirals are discontinuous with respect to $\mathcal{L}_{c}$
\end{diskontinuitet}

Given d U BV, we have a determined tangent, even when dU is  discontinuous (\cite{Riesz56}). A parameterization of spirals does not have a determined tangent, in particular when $(x,y)=(t x_{0}, t y_{0})$, we have $d U(x,y)=t d U(x_{0},y_{0})$, that is $\mid d U \mid=1$ implies t=1.

\subsubsection{Semi-regularity}

Assume $B(\phi,\psi)$ a bilinear form, then we have that the zero space $N=\{ (\phi,\psi)$ $ \quad B(\phi,\psi) = 0 \}$ defines orthogonals. Assume B absolute continuous, then we have d B=0 implies $B=const.$
B considered over a convex set, is  separately convex, that is separately absolute continuous. Consider for this reason T(N) and N completing over the spiral. Assume $(\phi,\psi)$ decomposable and p a semi-norm defined relative B.
Given $0 \in T(N)$ and p is continuous in 0, then p defines a norm on T(N). We assume a separation condition, for instance $\{ B < \lambda \} \subset \subset \Omega$. Given B separately locally 1-1 (with respect to norm), we have ${}^{t} B(\phi,\psi)=B(\psi,\phi)$ separately locally surjective. Note that, $B(\phi+\psi,\phi - \psi)=B(\phi,\phi) - B(\psi,\psi)$, given $B(\phi,\psi) \simeq B(\psi,\phi)$ locally. 

Assume $I = R_{1}R_{2}$, where  $R_{j}$ denotes reflection with respect to a bilinear form $B(f,g)$, for instance $R_{1}B(f,g)=B(-f,g)=-B(f,g)=BR_{1}(f,g)$. B can then be considered as semi algebraic, if for instance $BR_{j}(f,g)=R_{j}B(f,g)$, $j=1,2$ and $R_{j}^{2}=1$. Given BI=IB, then B is algebraic. Assume $< <f, d U> , g>= <f , < d {}^{t} U, g>>$, then we can define $<< f,g>, d U>$. Cf. the lifting principle, that is over analytic poly cylinders, for instance over $f,g \in L^{1}$, we have existence of U with $< f,d U>=0$ and existence of ${}^{t} U$ with
$< f, {}^{t} Ug>= <U f,g>$ and $<d {}^{t} U,g>=0$. Assume $U R_{1} = R_{1} U$ iff ${}^{t} U R_{2} = R_{2} {}^{t} U$, then we have $U {}^{t} U I = I U {}^{t} U$. Given R(U) normal and given the approximation property, we have $U$ is algebraic.

Assume $\int t(x,y) u(x) v(y) d x d y=0$, given $u \bot v$ relative $t d x d y$, that is $T(u \bot v)=0$ and given T of local type, we have T preserves orthogonality. 
\emph{Assume $T^{\bot}(u \otimes v)=\int t(x,y)u(x)v(y) d x d y=0$, given $t$ algebraic over $\mathcal{D}_{L^{1}} \times \mathcal{D}_{L^{1}}$, then $u \bot v$. That is,
$T(u \bot v)=T^{\bot}(u \otimes v)$. Assume $\mathcal{H}$ a normal space of distributions with
$\mathcal{H}'$ nuclear, then we have $T^{\bot}(u \otimes v)=0$ implies $u \bot v$ with respect to $t d x d y$.
Given t algebraic, the Lebesgue measure is absolute continuous relative $t d x d y$.} 

 Note (\cite{AhlforsSario60} Weyl's lemma) given $d U \bot T$ and $d U^{\bot} \bot \widehat{T}$ in $L^{2}$ and $T \bot \widehat{T}$, then we have $U \equiv U^{\bot}$ implies $d U(x) \in C^{\infty}$.
Assume $K(U,U^{\bot})=0$ a convolution kernel, locally 1-1, then we have $U=U^{\bot}$, iff ${}^{t}K$
surjective, that is we can identify the spiral. Consider for instance $C_{1} \parallel \phi \parallel^{\sigma_{1}} \leq \mid <K,\phi \otimes \psi> \mid \leq C_{2} \parallel \psi \parallel^{\sigma_{2}}$, where  $K$ is symmetric, that is K is invertible in $\phi,\psi$.

\emph{Assume $V=V_{1} \cup V_{2}$ with $V_{j}$ compacts in $\mathbf{R}^{n}$, then the convolution operator $\{ L \} : \mathcal{D}_{V_{1}}^{0}  \rightarrow \mathcal{D}_{V_{2}}^{m}$ is nuclear  iff 
$L \in C^{m}$} (\cite{Schwartz58}, Prop. 23). Note (\cite{Treves67}) that the polynomials are dense in $C^{m}$.
Consider in particular $V_{1}=\{ d U < \lambda \}$ and $V_{2}=\{ d U^{\bot} < \lambda' \}$, for real constants $\lambda,\lambda'$, that is we assume $U,U^{\bot}$ analytic in subsets $\Omega_{1} \subset V_{1}$ and $\Omega_{2} \subset V_{2}$ respectively. Assume further $\Omega_{1} \cap \Omega_{2} \neq \emptyset$.
Given $d U^{\bot}=\alpha d U$, with $\alpha \in \dot{B}$ and $U^{\bot}(f)= L* f$, we have $V_{1} \subset V_{2}$ and L can be approximated with polynomials. Note that when $U^{\bot}$ decreasing, compactness for $V_{2}$ is not immediate from the separation conditions.

According to prop 23 bis (\cite{Schwartz58}), without the condition $\Omega_{1}$ a subset of a domain of holomorphy, we have that we can do a nuclear continuation with loss of n+1 derivatives.

\subsubsection{Symmetry}
Assume $I_{x}(\varphi) - \gamma_{x} \in C^{\infty}$ and $I_{y}(\varphi) - \gamma_{y} \in C^{\infty}$, then we have $\varphi(x) - \varphi(y) \in C^{\infty}$. Given this implies $\varphi(x-y) \in C^{\infty}$, we have that separately regular with respect to x,y, implies very regular. In particular, consider $\varphi(x) - \varphi(y) \in \dot{B}$, when $x,y \rightarrow \infty$, implies $\varphi(x-y) \in \dot{B}$.

The symmetry operator ${}^{s} T(\varphi)=T({}^{s^{-1}} \varphi)$, does not imply $s^{2}=id$. Assume for instance S a multivalued domain with leafs $\{ S_{j} \}$, then we have that reflection on every leaf such that $s^{2}=id$, implies symmetry on S. Further, assume $U-I$ symmetric iff $U^{\bot}$ symmetric, where  S is generated by $U^{\bot}$.  
Assume $U$ absolute continuous, is defined by $X_{U}(f)=0$, then according to Hurwitz, we have that U is holomorphic or U=I, with respect to uniform convergence on compact sets. Given U projective $U^{\bot}=I-U$, we can define $(U,U^{\bot}) \underrightarrow{s} (U^{\bot},U)$. 

Assume $U^{\bot}=(I-U) + V$ and $\int (1 + \alpha) d U = \int \frac{d}{d T}(I + V)$. Given $I + V$ subnuclear, we have that convexity is preserved. The polar is given by V. Given every two points in R(V) can be combined through a path, a leaf is formed. 
In analogy with homology, we can write 
$d V^{\bot}=\Sigma \alpha_{j} d U_{2}$, where  the right hand side traces a polycylinder in the domain. We assume $d V^{\bot}$ absolute continuous relative $d U_{2}$.  Given $V \rightarrow I$ on the cylinder web, the limit is on one or several of the leafs. 

Assume $s=s_{1} s_{2}$, for instance $s_{1}$ harmonic conjugation and $s_{2}$ complex conjugation. Then we have $z \rightarrow \overline{z} \rightarrow \overline{z}^{*}$ gives $(x,y) \rightarrow (y,x)$ and $s^{2}=id$ that is reflection through the diagonal. When $s^{2}=s$, that is s a projection operator, we get the diagonal (the spiral when x,y is  $u_{1},u_{2})$. Starting from the two mirror model, where  $s_{j}$ is reflection through $L_{j}$, $j=1,2$, we have that F is  symmetric with respect to $L_{1},L_{2}$ and $L_{1} \rightarrow L_{2}$ projective, implies $F$ symmetric outside the diagonal. Assume $s=s_{1}s_{2}$, such that $s^{2}=s$, $s_{1}^{2}=-s_{2}^{2}=-I$, then we have $s_{2}s_{1}=-I$. Assume $s_{1}(A)=B$, $s_{0}$ the mapping over the axis $L_{1}$, then $s_{2}s_{0} s_{1}$ gives a normal model. When s is  projective in the plane, we can through $s_{0}$ compare with a transversal model. Necessary for the two mirror model is a separated space, that is points can be chosen outside the reflection axes.
 
Consider $y/x \rightarrow x/y$ regular, that is $y(x) \rightarrow x(y)$ invertible. Assume F subnuclear  and consider
$F_{1}(x,\frac{y}{x})$ with support on half the disk, in the same manner $F_{2}(y,\frac{x}{y})$ with support on the other half. Given $F=F_{1} + F_{2}$ and $I F_{1}=F_{2} I$,$I F_{2}=F_{1}I$, we have $I F=F I$.

\subsubsection{Sub-nuclear mappings} 
 
Assume $d \mid U \mid (\mid \varphi \mid)=0$ implies $\mid U \mid (\varphi)=\mid \varphi \mid$. Then $\mid \mid U \mid - \mid I \mid \mid \leq \mid U - I \mid$, that is U with constant variation does not imply U absolute continuous, even when $\mid U \mid$ absolute continuous. We assume that
$\mid U \varphi - \varphi \mid \rightarrow 0$ implies $U(\mid \varphi \mid) \rightarrow \mid \varphi \mid$. Note that $d \mid U \mid (\mid \varphi \mid)=0$
implies $\mid d U(\varphi) \mid =0$ implies U analytic over some $\varphi$.
Further, $\int_{\Omega} d \mid U \mid^{2}=0$, implies $\mid U \mid^{2} \mid_{\Omega}=0$, but $N(U)$ can have a positive measure (\cite{Riesz56})

Let $U^{\bot}=(I-U) + V$, then we have for instance $\widehat{U_{1} f} = U_{2} \widehat{f}= ((I - U_{1}) + V) \widehat{f}$. Over the lineality, we have $(I-U_{1}) \widehat{f}=0$ and V rotation.
Assume $\mathcal{H} \bigoplus \mathcal{H}_{1}=\mathcal{K}$ and $\Delta$ a restriction, with
$\Delta + (I-\Delta)=I$ and $R(I_{I-\Delta})=\mathcal{H}_{1}$. Assume $R(I_{\Delta}) \subset C^{\infty}$.
Given $I_{\Delta} + I - I_{\Delta}=I_{\Delta} + I_{\Delta - I}$,
assume $\Delta \in C^{\infty}(x,y)$ symmetric, with $\int \Delta(x,y) f(y) d y \equiv 0$, when $f \rightarrow \delta$, then we have $\Delta(0,0) =0$. Given $f \in \mathcal{D}$, we have for all derivatives, $I_{\Delta}(D^{\alpha} f) =0$ in $C^{\infty}$ and in the same manner for ${}^{t} \Delta$, why $\Delta \equiv 0$ (nbhd 0).

\subsubsection{Convexity}
Assume $\Omega=\{ d U(f)=0 \}$ and construct B, the absolute continuous closure of $\Omega$. Then we have that $B$ can be defined by $\{ U=I \}$. 
When $\Omega$ contains a 1-polar, we have that B has a 1-polar, but not conversely, that is if $U$ is  projective on $\Omega$, it is not necessarily projective on B.
Assume U 1-homogeneous and consider $z(t)=t f + (1-t) g$, then we have $U z(t)=t U f + (U g - t U g)$, thus convexity is preserved. But over the set where  $U t g=g$ (polar set), this is not the case. Thus, given U absolute continuous, we have that U is not necessarily subnuclear over continuations to g, such that $d U(g)=0$.
On B we have that $\int f' = f$ and given f=0 on $bd B$, we have $<f, d U>=\int {}^{t}( U f)' d x$. Thus, given ${}^{t} (U f)$ absolute continuous, we have $<f,d U>={}^{t} U f$. Given U subnuclear,  when $U \rightarrow I$, B is nuclear .

Assume $R(U)$ is defined by $d U$ and that $(d U)^{\bot}$ is given by $R(d U^{\bot})$. Assume the polar is given by dV, where $V \in \mathcal{G}$  and dV of distributional order 0. Assume $V=W_{1} W_{2}$, where $W_{j}$ are not spiral and $V^{\bot}=V$, then must we have $W_{1}^{\bot} \neq W_{1}$, why $W_{1}^{\bot}=W_{2}$. Note that
for the spiral to be well defined, the orientations for $W_{1},W_{2}$ must be compatible.  Given $U,U^{\bot}$ absolute continuous, we have $d (U - U^{\bot})(T)=0$ implies $U^{\bot}=(I-U)$, that is U projective over $R(T)$.

I is  optimally non-reduced, I preserves convex functions, given I preserves character over tangents, for instance $U \rightarrow I$ regularly, with $U f' = (Uf)'$.
 Assume L a connected line between f,g., ${}^{t} (U L)(f,g) ={}^{t} L {}^{t} U (f,g)$, that is sufficient for a connected line to be preserved under $U \rightarrow {}^{t} U$, is that $ {}^{t} L \rightarrow L$ is continuous.

\newtheorem{lyft}[max_order]{The Lifting principle}
\begin{lyft}
Assume $\Gamma$ chains of maximal order, without invariant sets and $\Gamma_{0}$ chains of order  0, that is cylindrical sets. On $\Gamma$ the composition of movements is not necessarily invertible. Through the lifting principle over $\Gamma_{0}$, we have existence of U, analytic over f, such that $g=Uf$.
\end{lyft}

Given $\Gamma_{0}$ has the approximation property and $I \in \overline{\mathcal{L}}_{c}(\Gamma,\Gamma_{0})$, then $\Gamma$ has the approximation property (\cite{Schwartz57}, prop. 2), that is $\Gamma_{0}$ has convex (analytic) continuation to $\Gamma$.

\subsection{Discontinuous continuation}
\subsubsection{Equilibrated sets}
For an equilibrated set, we have $e^{\lambda \varphi} \sim e^{\varphi}$, when $\mid \lambda \mid < 1$, that is $(I)(U) \ni e^{\lambda \varphi}=I(e^{\lambda \varphi})$
iff $\widehat{I}(\varphi^{\lambda}) =\varphi^{\lambda} \in (\widehat{I})$, that is we can assume this ideal is radical. Note that
Mackey topology does not imply absence of a spiral (over the polar), that is $(L_{c}')_{c}' \simeq L$ is  reflexive but does not imply projectivity.

Note that a closed contour can be compared with $z \rightarrow 1/z$ bijective. Assume $R_{j}$ $j=1,2$, such that $R_{1}R_{2}=I$ and $R_{1}U = U R_{2}$. 
Note $\mid I U \mid = \mid R_{1} U \mid = \mid R_{2} U \mid = \mid I U \mid$,
that is with respect to hypo continuity, $\mid U \mid$ is algebraic. Given $d \mid  U \mid = \alpha d \mid  U_{1} \mid$, with $\alpha \in L^{1}$, where  $\mid U_{1} \mid$ is analytic, $\mid U \mid$ is harmonic, that is has a decomposition into conjugated movements. 
Note that $\{ d \mid U \mid < \lambda \} \subset \{ d U < \lambda \} \subset \subset \Omega$, that is the separation condition implies hypo continuity, in the sense that $(\mid f \mid, d \mid U \mid(f))$ is  closed. 
Assume $ \mid F(\phi) \mid=G(\mid \phi \mid)$, where  G continuous over $\mid \phi \mid$, then F is $\sim$ hypo continuous. Given $F \rightarrow I$, we have $G \rightarrow I$. Assume $F=AB$ and dF absolute continuous 
with respect to d B, when $F \rightarrow I$, then A is invertible with respect to B. When I is subnuclear, we have existence of A with
$\frac{d A}{d B} \in L^{1}(d B)$, such that d F is absolute continuous with respect to dA.

When $u \otimes v$ is an absolute continuous function in $E^{'*}$, then we have given $d (u \otimes v)=0$, that $u(x)v(y)=1$, that is we can define $x \sim y$.
Given $u \otimes v(x,y)=\Phi(x-y)$ with $\Phi \equiv 0$ outside $x=y$, we must have $\Phi=\delta_{0}$.
Note that since $\mathcal{D}_{A} \times \mathcal{D}_{B} = \mathcal{D}_{A \times B}$ and when $u \otimes v(x,y)=\Phi(x-y)$, with $\Phi$ according to above, we can determine $x \sim y$ uniquely.

Assume $\Omega=\{ x \quad f(x) \in L^{1} \}$ and 
consider $\Omega_{ac}=\{ x \quad f(x) $ $\mbox{ absolute continuous } \} $ $\subset \Omega$ and $\Gamma_{ac}=\{ f \quad f \in L^{1}(\Omega_{ac}) \}$ . Given $f \in L^{1}$ there is $\Omega_{ac}(f)$, such that $f \in \Gamma_{ac}$. Over $\Omega_{ac}$ we have $\int f d I \sim \int f' d x$. Given a maximum-principle on $\Gamma_{ac}$ and $I \in \overline{\mathcal{L}}_{c}(G,\Gamma_{ac})$, we have that G has the approximation property.

Assume $T = P(D) d \mu$, where  $d \mu \in \mathcal{E}^{' (0)}$, then we have $T \in \mathcal{D}^{' m}$, given P of degree m.
Assume $\Gamma$ boundary for a pseudo convex domain, that is of order 0 with $(\Gamma)$ of finite order . This means that $(\Gamma) \subset \mathcal{D}^{' m}(E)$ can be generated by chains of finite order.
Given $d \mu$ of finite type, when we consider measures modulo $C^{\infty}$, $d \mu$ is of type 0. Further, $d \mu^{-1} \sim \frac{1}{1+V} d x \sim \Sigma c_{j} V^{j}$ is of finite order on $(\Gamma)$.
We assume invertibility for T only on $\Gamma$.

Given $d v=<f,e'> dx$, that is $\exists \frac{d v}{d x}=<f,e'>$, then through partial integration $\int \frac{d v}{d x} \varphi d x=\int <f,e'> \varphi(x) d x$.
Note $<f,dU> =<f, dU^{\bot}>$ does not imply $d U = d U^{\bot}$, when the space is not nuclear. Note that when $<f,d U_{1}><f,d U_{2}>(\varphi)=<f(\varphi),d U_{1} \otimes d U_{2}>=0$
implies $f (\varphi) \in C^{\infty}$, this does not give the same implication for movements separately. For instance, $(f \otimes g)^{\circ}=< X, f \otimes g >=0$ does not imply $X(f)=0$ or X(g)=0, for instance $X(f^{2}) \neq X(f)^{2}$.

Assume M of order 0 and P real, $S=P(D_{x}) M(x,y)$, then we have $S \circ S \sim P^{2} M \circ M$. Thus, $S^{N}$ is not locally of finite order, when $N \rightarrow \infty$, that is does not have support restricted to the diagonal. However $e^{P} -I=0$ iff $P=0$, that is when $S \sim e^{P}$ is  projective, $S^{\bot}$ must be locally of finite order.

\subsubsection{Separated spaces}

\newtheorem{hypo_algebra}[max_order]{Hypocontinuous convolution algebras}
\begin{hypo_algebra}
Assume $\mathcal{H}(\mathcal{G})$ a discontinuous convolution algebra, $U T \in \mathcal{H}$ implies $U^{-1} T \in \mathcal{H}$, $I T \in \mathcal{H}$ implies $(I-U) T \in \mathcal{H}$ but does not imply $(I - U)^{-1} T \in \mathcal{H}$. A hypocontinuous convolution algebra, implies that $(U,U^{\bot})$ has relatively compact sub level surfaces, given $U,U^{\bot} \in \mathcal{G}$ equicontinuous. 
\end{hypo_algebra}

More precisely, assume $\mathcal{H} \ni U T \rightarrow U^{\bot} T \in \mathcal{H}$, where  given $\mathcal{H}$ has the approximation property, we can define $U^{\bot}$ through a bilinear  form B. Assume p,q semi-norms over $\mathcal{H}$ relative $U,U^{\bot}$. Assume r a semi-norm relative $(U,U^{\bot})$. Then hypo continuity means (\cite{Treves67}), that $r(T(u,u^{\bot})) \leq p(T(u)) q(T(u^{\bot}))$. Assume $U^{\bot}=(I-U) + V$, then we have $d (I + V)=d U + d U^{\bot}=(1 + \beta) d U$. Given $(1 + \beta)^{-1} \rightarrow 0$, when $\mid x \mid,\mid y \mid \rightarrow \infty$, for instance $(1 + \beta)^{-1} \in \mathcal{D}_{L^{1}}$, then we have $\{ \beta < \lambda \} \subset \subset \Omega$, where  $\Omega$ an open set. Note if U is  algebraic over $\phi$ and $U \overline{\phi} = \overline{V \phi}$, we assume $U \sim V$, that is of the same character . With these conditions U is determined through dependence of the polar.

Assume p a semi-norm on $\mathcal{D}'$, q a semi-norm on $\mathcal{H}$, then $\mathcal{H} \subset \mathcal{D}'$ means that $p \leq q$. The property 
$(\epsilon)$ means that $p \sim q$. (\cite{Schwartz57}, Prop. 13).
 Assume $H_{p} \subset H_{w} \subset H_{q}$, where the weights are related to $L^{p}-$ spaces. Given p,q polynomials (algebraic),
it is not necessary that w is algebraic. However, given the property ($\epsilon$), w must be algebraic. 
Note that given $pq \leq p^{2} + q^{2}$, we can have a ``spiral''. 
In particular, precense  of spiral according to $(I) \subset (J)^{\bot}$ and $(J) \subset (I)^{\bot}$, when $(I)=R(U)$ and $(J)^{\bot}=R(U^{\bot})$, that is $(I)=(J)^{\bot}$.
Note $\mid <f,\widehat{g}> \mid \leq \parallel f \parallel \parallel \widehat{g} \parallel$,
that is given either $\parallel U f \parallel=0$ or $\parallel U^{\bot} \widehat{g} \parallel=0$, we have $U f \bot U^{\bot} \widehat{g}$, that is $R(U^{\bot}) \subset R(U)^{\bot}$.
Given $U f \neq 0$ and U locally 1-1, we have that $U^{\bot}$ surjective over $\widehat{g}$ and
$0 \neq U^{\bot} \widehat{g} \bot U f$, for some $\widehat{g}$.

Consider the completion to a closed curve, according to $\int_{\tilde{\gamma}} d c = \int_{\gamma} d \tilde{c}$, where  $\tilde{\gamma}$ a closed curve.
Assume $\gamma$ a spiral (not $\sim 0$), with $\tilde{c}=1$ over $\gamma$ (absolute continuous), given the spiral intersects a closed curve, we have that $\exists d c$ a closed form associated to $\tilde{\gamma}$. Every curve on the cylinder web, can be intersected by a spiral. But we do not have that $\tilde{\gamma} \sim 0$

\subsubsection{Uniformity}

When E is a barrel and $M' \subset E'$, we have $M'$ is weakly bounded iff $M' \subset \subset E'$ iff $M'$ equi-bounded (\cite{Schwartz57}, Prop. 2). Assume for instance $R(d U^{\bot}) \subset \subset R(d U)^{\bot}$. When $M' = \{ P < \lambda \}$ and $E'=\{ Q < \lambda \}$, the condition means $Q / P \rightarrow 0$ in $\infty$. The condition means that $\mid x \mid^{\sigma} \frac{Q}{P} < c$,
for some $\sigma > 0$, $\mid x \mid \rightarrow \infty$, that is the condition can be written $\{ P < Q \}$ relatively compact. 

Assume $d \mu$ BV on L completed and $M=\{ \varphi \in L \quad \int \mid d \mu \mid(\varphi)=const \}$. Assume $d V=\mid d \mu \mid$ and $d V$ absolute continuous relative d I on M, then we have $V$ absolute continuous on M. That is, $d V= \alpha d I$, where $\alpha \in L^{1}_{ac}$ on M and $V= const. I$ on M. 

Assume $d F \subset d U \otimes d U^{\bot}$, with $(d U)^{\bot} \subset (d U^{\bot})$. Given $F$ absolute continuous, we have projectivity over analytic chains. Given $dU  \bigoplus U^{\bot} \sim 0$, implies $U + U^{\bot} = I$ outside the polar, where for instance $I (U \rightarrow U^{\bot}) \notin \mathcal{L}_{ac}$. Over absolute continuous functions, we have thus $d U \sim d I$, that is dU approximates d I regularly. Note that precense of a regular chain is  necessary
for reduction to absolute continuous functions.

Given a barrel T, there is a semi-norm p, such that the unit ball can be given by $\{ p \leq 1 \}$. p is  continuous iff $0 \in T$ (\cite{Treves67}).
Thus, the semi-norm corresponding to the cylinder web is not continuous. The approximation property through regularizing and truncation on a normal space implies the strict approximation property (\cite{Schwartz57}, Prop. 3, pg. 10).

Assume $W_{c}=\{ (x,y) \quad \mid d U(x,y) \mid = c \}$ compact and $(W_{c})=\{ (x,y) \quad \mid d U(x,y) \mid < c \}$. Given $(W_{c}) \subset \subset \Omega$, with $\infty \in \Omega$, then $\mid d U \mid$ can be given as locally reduced. When $\mid d U \mid=\beta d I$, where $\beta$ is a continuous function, maximum is assumed on $W_{c}$. Note that $d U$ is not necessarily locally reduced. 

\emph{Assume $<<T,\varphi>,e'>=<T,e'>(\varphi)=0$ for arbitrary $\varphi \in \mathcal{D}_{L^{1}}$, then using the approximation property for $\mathcal{D}_{L^{1}}'$ and relative topology for uniform convergence on compact sets, we have existence of $d F$ closed locally, with $\int <T,e'> \varphi d x =  \int \varphi d F$. In particular, when $e'$ defines an analytic movement, $\int <T,e'>$ is on closed contours locally independent of choice of local coordinates.}

\subsubsection{Discontinuous continuations}
Consider $(U,U^{\bot})$ as a continuation element, that is $R(U)=R(U,0) \subset R(U,U^{\bot})$.
Given $(U,U^{\bot})$ locally 1-1 , we have given Schwartz separation condition, ${}^{t} (U,U^{\bot}) \simeq (U,U^{\bot})^{\bot} \simeq (U^{\bot},U^{\bot \bot})$. Given U reflexive, then $(U,U^{\bot}) \simeq U + U^{\bot}$ can be considered as projective. Given $U^{\bot}=(I-U) + V$ and U absolute continuous, we have $d U=0$ implies $U^{\bot}=V$ on the polar. Further, when $d U^{\bot}=0$, we have $V^{\bot}=V$ on the polar. 

Assume $u \rightarrow v$ decomposable, in the sense that we have corresponding semi-norms $p_{u} \geq q_{v}$. 
Given $d \mid U^{\bot} \mid = \beta d \mid U \mid$ and $d U^{\bot}=\alpha d U$, we have that $\beta=const$ does not imply $\alpha=const$. Given $V(U,U^{\bot})$ the polar, we have that $R(V(U,U^{\bot}))$ is decomposable does not imply $R(V(\mid U \mid, \mid U^{\bot}) \mid)$ decomposable. 

 \emph{Assume $K_{1},K_{2} \in \mathcal{H}$. Consider $K_{1}(l,u) K_{2}(m,v) \simeq K(l,m) K_{0}(u,v)$, where  $K_{0}$ regulates the change of variables $u \rightarrow v$, for instance ${}^{t }K_{2}(v,m) = v^{*} \otimes in {}^{t} K_{2}(u,m)$ (\cite{Schwartz58}). The condition $\frac{d v}{d u} \in L^{1}(d u)$ gives an integrable change of variables. Given $\mathcal{H}$ has the approximation property through truncation, we can define a continuous deformation.} 

In particular for instance $K_{1}(u,u) K_{2}(u^{*},w) \simeq K(u,u^{*})K_{0}(u,w)$, where $u,w$ are in complementary sets.
Assume $f,\widehat{g}$ analytical and $f \widehat{g} \equiv 0$ in nbhd $u_{0}$, then $N(f \widehat{g})=N(f) \cup N(\widehat{g})$. When $\{ f \widehat{g} < \lambda \}$ is semi-algebraic with $\widehat{g}(u)=h(u^{*})$ and when the set is decomposable, we can assume $f \widehat{g} \sim f(u) P(1/u)$,
for some polynomial P. More generally, assume $d I$ BV considered over ac functions, then we have $d U + d V = d I$, where $d U,d V$ are BV measures, $d U$ non-decreasing and $d V$ non-increasing. We can extend the domain to $\alpha d I$, where $\alpha \in L^{1}(du)$, to
motivate an integrable change $u^{*} \rightarrow w$.

Given $(I_{1}) \subset (I_{2})$, ideals of analytic functions, we have $N(I_{2}) \subset N(I_{1})$, but first surfaces to $I_{1}$ are not necessarily first surfaces to $I_{2}$.
For instance $(I_{1})$ has the strict approximation property, but $(I_{2})$ has only the approximation property. Assume $(I_{1}) \subset (I_{2})$ according to $p/q \rightarrow 0$ in $\infty$ for the corresponding semi-norms. Given $T_{j} \in (I_{1})$ implies $T_{j} \in (I_{2})$, that is $T_{j} \rightarrow T$ in $(I_{1})$, does not imply $T_{j} \rightarrow T$ in $(I_{2})$, for instance, $T_{j} \rightarrow T + c$ in $(I_{2})$.

$L^{1} \rightarrow L^{2}$ is not bijective, $d U^{2}$ locally 1-1 (absolute continuous) does not imply $dU$ locally 1-1 (absolute continuous). 
According to Lie (\cite{Lie91} chapter 12), we have $d U_{S}(f)=(y + \kappa x) \frac{\delta f}{\delta x} + (-x + \kappa y) \frac{\delta f}{\delta y}$. Assume $d U(f)= y \frac{\delta f}{\delta x} - x \frac{\delta f}{\delta y}$.
\emph{Thus,  $d U - d U_{S} \simeq - \kappa d U^{\diamondsuit}$ (harmonic conjugate).}
Assume $I \prec U=A+B$, with A invertible, then we have $A^{-1}(I + B/A)^{-1}=(A + B)^{-1}$. Necessary for convergence is that $\mid B/A \mid < 1$, but globally further $B/A \rightarrow 0$ in $\infty$, that is $B \prec \prec A$. 
$U_{S}=A+B$, with $A \sim B$ gives that \emph{the inverse to $U_{S}$ does not converge in an analytic sense of measures.}

\subsection{Spiral completion}

\subsubsection{The approximation property}
Assume $F$ has the approximation property and that $I \in \overline{\mathcal{L}_{c}}(E,F)$, implies that E has the approximation property. Thus $E \underrightarrow{U} F \underrightarrow{I} F$ and $E \underrightarrow{I} E \underrightarrow{U} F$, that is U is algebraic in $\mathcal{L}_{c}$.
For instance, assume $U^{\bot} = \{ U_{j} \}^{\bot}$ multivalent, and that
$d U_{j}=\alpha_{j} d U_{1}$ and $d U_{1}=\Sigma d U_{j} / \alpha_{j}$, where  $\alpha_{j} \rightarrow 0$ in $\infty$. If for instance the polar is given by $d V$ and $d U_{j} = \beta_{j} d V$, the condition for invertibility $d U_{1} = \gamma d V$, is $\gamma=\Sigma \beta_{j}/\alpha_{j} \rightarrow 0$ in $\infty$.

Assume $\mathcal{H}=\mathcal{H}_{0} + \mathcal{H}_{p}$, given $I \in \overline{\mathcal{L}_{c}}(\mathcal{H},\mathcal{H}_{0})$, where  $\mathcal{H}_{0}$ has the approximation property, we have the same for $\mathcal{H}$ (\cite{Schwartz57}). Note $I \in \overline{\mathcal{L}_{c}}$ does not imply that I is locally reduced. Assume $\mathcal{H}_{0}$ is defined by $(U,U^{\bot})$, with $U \neq U^{\bot}$. Note $U=U^{\bot}$ defines the web, convex with respect to spirals. 
The condition $F(U,U^{\bot})=0$ implies $U=U^{\bot}$, that is F is reduced over $(U,U^{\bot})$, can be written $\int \phi d F=0$ for $\phi \in (I)$ implies $(I) \subset C$ ($\phi \neq 0$).   Given $d F=0$ implies $F=I$, that is a constant surface, means that $U^{\bot}$ is not uniquely determined outside the polar.

Assume $d U = \alpha d V$, with $\alpha \rightarrow 0$ in $\infty$, for instance $\alpha \in \dot{B}$. We then have inclusion between the corresponding ideals, further
$\mathcal{L}_{c}(d V) \rightarrow \mathcal{L}_{c}(d U)$ is continuous, that is we have topology induced by truncation. Thus, existence of an uniformly convergent subsequence in $d V$ implies existence of an uniformly convergent subsequence in $d U$.  

Assume $U^{\bot}=U + V$. Given $\{ V=0 \}$ removable, we have regular approximations. Assume V is defined through completion to $L^{1}$, in the sense that dx is absolute continuous relative d V. Given $d V=0$ implies $-V=I$, then we have over f absolute continuous, that the polar is removable.
For a normal space and the closure of $U^{\bot}$ in $L^{1}$, we have that it is sufficient for a removable continuation, that $U^{\bot} \neq U$ close to $U=I$.

Consider hypo density, that is $\mid P \mid^{2}$ polynomial does not imply P polynomial. If T is downward bounded by $P_{1}P_{2}$, we do not necessarily have separate density. Consider for instance y=y(x), why $P_{1}(x)P_{2}(y)$
polynomial in $(x,y)$ does not imply $P_{1}P_{2}$ polynomial in x.

\subsubsection{Completing sets}
Given $\mathcal{K}$ normal with the approximation property through regularization, truncation, given $\mathcal{T}$ saturated and completing, with
$S={}^{t} \Delta T$ and ${}^{t} \Delta : \mathcal{H}' \rightarrow \mathcal{K}'$, then the scalar product is independent of T (\cite{Schwartz58}, Prop 20 ), for instance given $\Delta$ Fourier
dual, then $<\varphi,T>$ can be given through only $N_{\alpha}$ (\cite{Nilsson72}), $Exp_{\parallel \parallel_{L^{1}}}$ (\cite{Martineau}), when the movement is defined through action in phase (\cite{Lie96})

Given the conditions in prop. 20 (\cite{Schwartz58}), we can reduce the condition for hypoellipticity, from the plane of $\mid f \mid,\mid \widehat{f} \mid$ to plane of $f$, for instance it is necessary that f is algebraic
over the spiral (polar). 
\newtheorem{completing}[max_order]{Necessary condition for hypoellipticity}

\begin{completing}
Without the condition on the approximation property, it is necessary in the condition for hypoellipticity over maximal chains, to consider the spiral to a completing set
\end{completing}
(\cite{jag20}).
Given $T \in \mathcal{D}(F,\beta_{0})$ has support in A, then we have $(\mathcal{D}_{A})' \simeq \mathcal{D}' / (\mathcal{D}_{A})^{0}$, where  $\mathcal{D}_{A}^{0}$ is orthogonal to $\mathcal{D}_{A}$ in $\mathcal{D}'$. When $U^{\bot}$ is defined starting from the scalar  product, it is sufficient to define $U^{\bot}$ over $\widehat{\varphi}$, the product is not dependent of $\mathcal{F},\mathcal{H},\mathcal{K}$. Given $U^{0}$ the polar to U, we have $\mathcal{K}_{U^{0}}' \simeq (\widehat{\mathcal{K}}_{U})'$ (\cite{Schwartz58}).

Hilbert spaces have the metric approximation property. When $F_{B}$ a Banach-space (B completing), we have given the parallellogram law, that $F_{B}$ is Hilbert. When $\parallel d U - d U^{\bot} \parallel=0$ implies a measures zero set, then $d U$ is projective and $d U^{\bot}$ generates $(d U)^{\bot}$. Consider $\Gamma=\{ d U=d U^{\bot} \}$, given d U locally reduced, then $(d U)^{\bot}$ does not leave space for a spiral. Let $\tilde{U}=(U,U^{\bot})$.

\emph{Note if $R(\tilde{U})$ completing set, then we have through change of variables to $(U_{1},U_{2})$, that the radius L can be seen as representing a spiral, as long as the change of variables preserves decomposability.}

We consider $d (U,U^{\bot})=(\alpha,\beta) d (U_{1},U_{2})$, envelop to $R(\tilde{U})$ generated by $(U_{1},U_{2})$. Assume p,q semi-norms relative $U,U^{\bot}$, Consider r a semi-norm
relative the polar $V=U^{\bot} + U - I$. Assume $r(V) \leq p(U) q(U^{\bot})$. Over a convex set of f, given r defines a norm, we have that $r(V(f))=0$ implies V(f)=0, that is projectivity for $(U,U^{\bot})$.
Note that when the domain is defined by $(U_{1},U_{2})$, we have for the polar V that $I + V \simeq U_{1} + U_{2} = U_{2} + U_{1}= I + V^{\bot}$, where  we assume reflexivity over the domain. 

Given existence of $\limsup \mid U \mid = +\infty$, we do not have necessarily existence of $\limsup U$ (for instance essential singularities). Sufficient,  to guarantee two-sided limits, is  that U is acting algebraically in phase (preserves constant value in $\infty$) and  $U(\frac{1}{\phi})=U^{-1}(\phi) \in L^{1}$. Note that $I=I_{red}$, implies $I(f) I(1/f) = I(1)=I$. $\gamma$ can be seen as a closed contour, given $z \rightarrow 1/z \rightarrow z$ injective, that is $\psi(z)=1/z$ and $\psi^{2}=I$.  
Assume $d U = \alpha d U_{1}$ with $\alpha \in L^{1}(dU_{1})$. Let $\Sigma_{1}=\{ d U_{1}=0 \}$
and $\Sigma=\{ d U=0 \}$. Further, $\Sigma_{1}^{0}=\{ \Delta U_{1}=0 \}$ and $\Sigma^{0}=\{ \Delta U =0 \}$ (Laplace operator). Associated to $\Sigma_{1}^{0}$, we have an conjugated resolution of identity. Assume
$\Omega=\{ \alpha=const \}$ and the closure $\tilde{\Omega}=\{ \alpha \in L_{ac}^{1} \}$. Note that a resolution starting from $\Omega$ has eigenvectors $\subset \subset \mathcal{H}$ (through linear  independence). Given $\tilde{\Omega} \cap \Sigma = \emptyset$, we have $\tilde{\Omega} \cap \tilde{\Sigma}^{0} \neq \emptyset$, where $\tilde{\Sigma}^{0}$ is the envelope to $\Sigma^{0}$ (\cite{jag19}). Given (nbhd $\tilde{\Omega}$  - preimage P=0, for a polynomial P) convex, we can relate to a regular  covering. 
Assume $d U^{\bot}=\beta d U^{\diamondsuit}$ (harmonic conjugate) and that the mapping corresponding
to $\Sigma^{\diamondsuit} \rightarrow \Sigma^{\bot}$, maps lines on lines 1-1 and maps closed forms on closed forms 1-1. Thus, where  $d U, d U^{\bot}$ are analytic, $d U$ is harmonic.

\subsubsection{The spiral}
An approximation with sequential movements of a spiral, has a discrete intersection with the spiral. The cylinder web is  convex for spirals, every sequential movement in $(U_{1},U_{2})$ is  in the complement to the spiral. The spiral is  a purely complex figure in $(u_{1},u_{2})$, transversals purely real, that is $F(u_{1},u_{2})$ on spirals is  complete. Spirals give precense of a complete approximation, in the case of holomorphic leafs (\cite{Oka60} Puiseux approximation).

Assume $\mathcal{G}=\mathcal{G}_{1} \bigoplus \ldots \bigoplus \mathcal{G}_{n}$. Given $\Omega$ a domain, such that $U_{j}$ is analytic over $\Omega$ for all j, we have a covering. This is  not evident for a very regular boundary, where only some $U_{j}$ is analytic. 

Assume $\tilde{I}(\Omega)=I(T(\Omega))$. Then we have $\int f d U \simeq \int \frac{d}{d t}({}^{t} U f) d t$ and given ${}^{t}Uf$ is absolute continuous, $\simeq {}^{t} U f$ (boundary condition ${}^{t} U f=0$ on bd T ( $\Omega$)).
Assume U analytic over f, that is log U f convex. Assume $\tilde{I} \simeq \mbox{ rad }I$, then we have that $\phi_{N}' \simeq \phi_{1}'$, where  $\phi_{N}$ is the phase corresponding to $f^{N}$ and $\phi_{N}$ is convex iff $\phi_{1}$ convex.

Convexity with respect to path: Assume $ d U=\alpha d V$, where $\alpha \in \dot{B}$, then we have that U is convex with respect to V, given $\frac{d \alpha}{d V} \geq 0$. Given $L(f,g)=t f + (1-t) g$, we have that $\frac{d L}{d t}=f-g=\int_{f}^{g} d \mu$. Given $\frac{d L}{d t}=0$, we have that L defines a closed path. Define a continuous deformation, $U f + U^{\bot} \widehat{g}=W_{T}(f,g)$,
where  $U f \rightarrow f$ iff $U^{\bot} \widehat{g}=0$, then
given $f,g \in R(W)$ locally convex, we have $W_{T} (f,g) \subset R(W)$ and we have precense  av intermediate values $h \in R(W)$ under the condition $R(U)^{\bot} \subset R(U^{\bot})$. 

Starting from $U^{\bot}=(I-U) + V$, with U reflexive but not projective, we have that when $U \rightarrow I$, $U^{\bot}$ approximates a spiral and $V=V^{\bot}$ approximates $U^{\bot}=I$. The mapping $U \rightarrow U^{\bot}$ assumes regularity, if $Uf \sim e^{\phi}$, then $\{ \phi = - \infty \}$ is mapped on the set $\{ \phi=0 \}$, that is (U - I)f=0. Given $\phi$ locally algebraic, the polar must be the preimage of a locally algebraic set.

Given E is separable, then T is locally summable iff we have existence of $f \in E^{' *}$, such that $U T(\varphi)=\int <f(x),d {}^{t} U> \varphi (x) d x$. (\cite{Riesz56}, Radon -Nikodym).
Consider $w T$ locally bounded in $\mathcal{H}$, where w is a multiplier or weight. When w is reduced, we have that T is locally bounded. Alternatively, $w \widehat{T} \rightarrow 0$, defines T relative the $L^{1}$-norm.
Note that $\{ \alpha=0 \}$ has normals of the type $\{ \alpha = const \}$, given the separations condition (isolated singularities). 
For instance, assume S a measure on $\Gamma$, such that $S-I$ is reduced modulo $\Gamma$, with $S \in \mathcal{E}^{'(0)}((\Gamma))$, that is $(S-I)(\phi)=0$ implies $\phi \in \Gamma$. Given S reduced and convex (absolute continuous) on $\Gamma$, it is sufficient that d S=0 on $\Gamma$. 
Note that the cylinder web C is convex with respect to $d U_{S}$ and $d U_{S} \varphi \rightarrow \varphi$
implies $\varphi \in C$,  but $d U_{S} \nrightarrow 0$, that is $0 \notin C$. However $0 \in T(C)$, the convex disk closure. When $d U_{S}$ is a measure on $C$ and $\tilde{d U_{S}}$ a subnuclear measure on $T(C)$, we have that $\tilde{d U}_{S} \varphi \rightarrow \varphi $ implies $\varphi \in T(C)$. Note $d U_{S}=\xi \delta_{x} + \eta \delta_{y}$, has $\xi,\eta \rightarrow 0$, when $(x,y) \rightarrow 0$.

\bibliographystyle{amsplain}
\bibliography{vektor2021}

\end{document}